\documentclass[12pt,a4paper]{article}
\usepackage[centering]{geometry}
\usepackage{blindtext}
\usepackage{amsmath}
\usepackage{setspace}
\usepackage[center]{titlesec}
\usepackage{caption}
\usepackage{multirow}
\usepackage{graphicx,xcolor}
\usepackage{subfigure}
\usepackage{bm}
\usepackage{cite}
\usepackage{amsfonts}
\usepackage{float}
\usepackage{mathrsfs}
\usepackage{listings}
\usepackage{verbatim}
\usepackage{abstract}
\usepackage{amsthm}
\usepackage{amssymb}
\usepackage{lastpage}
\usepackage{fancyhdr}
\usepackage{epsfig}
\usepackage{epstopdf}
\usepackage{sectsty}
\usepackage{titlesec}
\usepackage{authblk}
\usepackage{indentfirst}
\titleformat{\title}{\centering\bfseries\fontsize{16pt}}
\titleformat{\section}{\Large\bfseries\flushleft}
\titleformat{\subsection}{\Large\bfseries\flushleft}
\titleformat{\subsubsection}{\Large\bfseries\flushleft}

\theoremstyle{plain}

\theoremstyle{definition}

\theoremstyle{remark}
\newtheorem{rem}{Remark}

\providecommand{\keywords}[1]{\textbf{Keywords:} #1}

\title{Mixed GMsFEM for linear poroelasticity problems in heterogeneous porous media}
\author{Xia Wang\thanks{E-mail address: xwang@math.cuhk.edu.hk. Department of Mathematics, The Chinese University of Hong Kong, Shatin, Hong Kong SAR.}, Eric Chung\thanks{E-mail address: tschung@math.cuhk.edu.hk. Department of Mathematics, The Chinese University of Hong Kong, Shatin, Hong Kong SAR.}, Shubin Fu\thanks{E-mail address: shubinfu89@gmail.com. Corresponding author. Department of Mathematics, University of Wisconsin-Madison, USA},
and Zhaoqin Huang\thanks{E-mail address: huangzhqin@upc.edu.cn. School of Petroleum Engineering
China University of Petroleum (East China), China.}
}

\begin{document}
\maketitle{}
\begin{abstract}
Accurate numerical simulations of interaction between fluid and solid play an important role in applications. The task is challenging in practical scenarios as the media are usually highly heterogeneous
with very large contrast. To overcome this computational challenge, various multiscale methods are developed. 
In this paper, we consider a class of linear poroelasticity problems in high contrast heterogeneous porous media, and 
develop a mixed generalized multiscale finite element method (GMsFEM) to obtain a fast computational method. 
Our aim is to develop a multiscale method that is robust with respect to the heterogeneities and contrast of the media, and gives a mass conservative fluid velocity field. 
We will construct decoupled multiscale basis functions for the elastic displacement as well as fluid velocity. Our multiscale basis functions are local.
The construction is based on some suitable choices of local snapshot spaces and local spectral decomposition, with the goal of extracting dominant modes of the solutions.
For the pressure, we will use piecewise constant approximation.
We will present several numerical examples to illustrate the performance of our method. 
Our results indicate that the proposed method is able to give accurate numerical solutions with a small degree of freedoms. 
%We utilize mixed finite element method, as well as following the framework of Generalized Multiscale Finite Element Method (GMsFEM). Two splitting schemes are studied for temporal approximation: fixed-stress and fully coupled. For the approximation of velocity field, we construct snapshot space and select the dominant basis functions by performing spectral problems. Multiscale space for displacement are established by implementing spectral decomposition directly on the local fine-scale basis functions. The basis functions obtained for velocity is mass-conservative. The waiver of snapshot space construction for displacement simplify the implementation procedure. Numerical experiments shows that our method can achieve good accuracy with only a few multiscale basis functions selected.

\end{abstract}
\keywords{multiscale method, mass conservation, poroelasticity}

\section{Introduction}
 
Simulation of interaction between fluid and solid constituents within a heterogeneous porous medium is of vital importance in areas such as reservoir geomechanics \cite{Gambolati2006, Zoback2010,deng2017locally} and medical diagnosis\cite{Mura2016}. The mechanical behavior of such porous media accounts for the coupling of the solid deformation and fluid flow behavior. Among all the models proposed by pioneering researchers, Biot introduced a three-dimensional theory of elastic deformation of fluid infiltrated media \cite{biot1941general} and extended it to porous media in 1956 \cite{Biot1956}, which can accurately model the dynamic behavior within porous media. 

Due to the presence of heterogeneity, direct simulation of the model problem requires a high grid resolution which is computationally expensive. There are in literature a number of multiscale methods for solving these problems with a reduced computational cost. Some popular examples are upscaling or homogenization approaches (e.g., \cite{durlofsky1991numerical, efendiev2009multiscale, wu2002analysis, gao2015numerical,owhadi2007metric}), heterogeneous multiscale methods (HMM) \cite{abdulle2012heterogeneous,weinan2003heterognous}, multiscale finite element methods (MsFEM) \cite{aarnes2004use, arbogast2007multiscale, efendiev2009multiscale, hou1997multiscale}, generalized multiscale finite element method (GMsFEM) (e.g.,\cite{efendiev2013generalized,chung2016adaptive,chung2016mixedperforated,chung2017onlineperforated,chan2016adaptive,yang2020online,yang2019multiscale}) and local orthogonal decomposition method (LOD) \cite{brown2016multiscale}. The goal of these approaches is to construct low dimensional computational models which can give approximate solutions with good accuracy. For instance, numerical homogenization aims at computing an effective quantity for the heterogeneous coefficient so that the resulting computational model can be solved on a coarse grid to give an upscaled solution. 
%which can preserve some averages for a given set of local boundary conditions and then solve the problem by a reduced dimension space on the coarse grid. 
%Homogenization method, however, can only apply to models that can be reduced to a low dimensional structure. In some cases, the limited effective dimension of microscale problems is not possible. 
Another way is to represent the solution by some carefully designed local multiscale basis functions as in MsFEM.
These basis functions are solutions of local problems with appropriate boundary conditions. Contrary to standard finite element basis, MsFEM basis are oscillatory in the interior of each coarse block, and these features are important in capturing oscillations in the solutions. Therefore, MsFEM basis functions contain more information and are good representatives of the solution space. Nevertheless, the accuracy of MsFEM depends on local boundary conditions and assumes scale separation. Though effective in many cases, multiscale methods that only use local information may not accurately capture the local features of the solution. GMsFEM is a generalization of MsFEM with the goal of designing a systematic way to enrich the multiscale solution space. 
It consists of two stages: offline stage and online stage. In the offline stage, we construct a small dimensional multiscale basis functions that can be effectively used to solve the global problem in the online stage for any input parameter, such as right-hand sides or boundary conditions. To get these small dimensional multiscale basis functions, we first compute snapshot spaces locally and then  reduce the snapshot space by performing a suitable spectral decomposition. The spectral problems are designed by error analysis and have a huge impact on the convergence rate of the method. In the online stage, basis functions can also be constructed and added based on the solution residual with aims of reducing error significantly and capturing global information \cite{yang2020online}. 

Our work is motivated by the framework of GMsFEM. There are in literature research on GMsFEM for poroelasticity problems (see e.g., \cite{brown2016generalized, fu2019computational}). The aim of our work is to handle the critical need of mass-conservation in flow problems. Several mixed finite element methods have been developed to cope with this challenge (see e.g., \cite{arbogast2007multiscale, chen2003mixed, chung2015mixed,yang2020online}). In multiscale framework, some mixed methods enjoy good property of mass-conservation  without post-processing. In consideration that our model is based on viscous flows and governed by Darcy's law, we introduce the velocity variable. Therefore, we are aiming to find appropriate space for displacement, velocity and pressure. For the approximation of velocity field, we first construct snapshot spaces which are local solutions supported on single coarse edge neighborhood and are consisting of all possible boundary condition of unit flux with respect to the fine grid. The offline space of velocity field is achieved by performing local spectral problems in the corresponding snapshot space. In the framework of continuous Galerkin approach, one basis function per edge is not sufficient to capture many disconnected multiscale features \cite{efendiev2013generalized, chung2014adaptive}, while our method can systematically generate enough basis functions to represent the multiscale features. Moreover, there is no need to use partition of unity functions. For approximation of displacement field, we use local fine basis functions as the snapshot functions. Spectral problems are performed to get the multiscale basis functions for displacement. For pressure basis functions, piecewise constant functions are proved to be good approximation in our numerical experiments. For time sequential approximation, we consider two splitting schemes as discussed in paper \cite{kolesov2014splitting, kim2009stability, ferronato2010fully}: fixed-stress and fully coupled. Fully coupled scheme generates a bigger matrix, while fixed-stress scheme is more economical.

The paper is organized as follows. In Section \ref{Sec:Preliminaries}, we introduce the poroelasticity model. We define the mesh and partition, derive the variational formulation, and apply different splitting schemes in Section \ref{Sec:variational and splitting}. Construction of multiscale velocity basis and multiscale displacement basis are presented in Section \ref{Sec:multiscalebasis}. In Section \ref{Sec:numerical results}, numerical results are illustrated, and we observe that our proposed method is able to give accurate solutions with a small dimensional approximation space. The paper ends with a conclusion.

\section{Preliminaries}\label{Sec:Preliminaries}
We let $\Omega\subset \mathbb{R}^d$ $(d=2,3)$ be a bounded computational domain with 
Lipschitz boundary. Let $T>0$ be a fixed time. 
We consider the following linear poroelasticity problem in which we find the displacement $u$ and the pressure $p$ satisfying
\begin{subequations}
\begin{align}
-\nabla\cdot \sigma(u)+\alpha \nabla p = 0 &\quad\text{in } (0,T]\times \Omega,\label{Pro:formula1a}\\
\alpha \frac{\partial \nabla\cdot u}{\partial t}+\frac{1}{M}\frac{\partial p}{\partial t}-\nabla\cdot\Big(\frac{\kappa}{\nu}\nabla p\Big)  = f &\quad\text{in }(0,T]\times\Omega,\label{Pro:formula1b}
\end{align}	
\label{Pro:formula1}
\end{subequations}with the initial condition $p|_{t=0} = p_0$ for the pressure.  We split the boundary of the domain into two parts $\partial \Omega = \Gamma_1\cup\Gamma_2$. We assume the following boundary conditions on each portion
$$u = 0, p=0 \text{ on }\, (0,T]\times \Gamma_1,\quad u = 0, \frac{\kappa}{\nu}\nabla p\cdot \vec{n} = 0 \text{ on } \, (0,T]\times \Gamma_2,$$
where $\vec{n}$ is the outward unit normal vector on $\partial\Omega$.
In Problem \eqref{Pro:formula1}, we denote the stress tensor by $\sigma(u)$, the Biot modulus $M$, the fluid viscosity $\nu$, the source term $f$, and  the Biot-Willis fluid-solid coupling coefficient $\alpha$. For models derived from linear elastic stress-strain constitutive relation, the stress tensor is expressed as 
$$\sigma(u) = 2\mu \epsilon(u)+\lambda \nabla\cdot(u)\mathcal{I},\,\, \epsilon(u)=\frac{1}{2}(\nabla u+\nabla u^T),$$
where $\mathcal{I}$ is the identity matrix, $\lambda,\mu>0$ are the Lam$\acute{e}$ coefficients. The Lam$\acute{e}$ coefficients can be expressed in terms of the Young's modulus $E>0$ and the Possion's ratio $\eta\in (-1,\frac{1}{2})$ via,
\begin{equation}
\lambda = \frac{\eta}{(1+\eta)(1-2\eta)}E,\quad \mu = \frac{1}{2(1+\eta)}E.
\label{eqn:relationlambda}
\end{equation}
Here the primary sources of the heterogeneities in the physical properties arise from $M, \lambda, \mu, \alpha$ and $\kappa$.
% arise from $M$ the Biot modulus,  $\lambda,\mu$ the Lam$\acute{e}$ coefficients,  $\alpha$ the Biot-Willis fluid-solid coupling coefficient and $\kappa$ the permeability.

To proceed with the mixed finite element method, we introduce the velocity variable $$g = -\frac{\kappa}{\nu}\nabla p$$ to Problem \eqref{Pro:formula1}.  To state it more clearly, we are dealing with the following  problem: find $(u,g,p)$ such that
\begin{subequations}
	\begin{align}
	-\nabla\cdot \sigma(u)+\alpha \nabla p &= 0\quad\text{in }(0,T]\times\Omega,\label{Pro:formula2a}\\
	\kappa^{-1}\nu g+\nabla p &= 0\quad\text{in }(0,T]\times\Omega,\label{Pro:formula2b}\\
	\alpha \frac{\partial \nabla\cdot u}{\partial t}+\frac{1}{M}\frac{\partial p}{\partial t}+\nabla\cdot g & = f\quad\text{in } (0,T]\times\Omega,\label{Pro:formula2c}
	\end{align}	
	\label{Pro:formula2}
\end{subequations}
with initial and boundary conditions rewritten as 
$$u = 0, p=0 \text{ on }\, (0,T]\times \Gamma_1,\quad u = 0, g\cdot \vec{n} = 0 \text{ on } \, (0,T]\times \Gamma_2.$$

\section{Variational Formulation and Splitting Scheme}\label{Sec:variational and splitting}
\begin{figure}[H]
	\centering
	\subfigure[Coarse grid and fine grid.]{
		\begin{minipage}[t]{0.4\linewidth}
			\centering
			\includegraphics[width=2.4in]{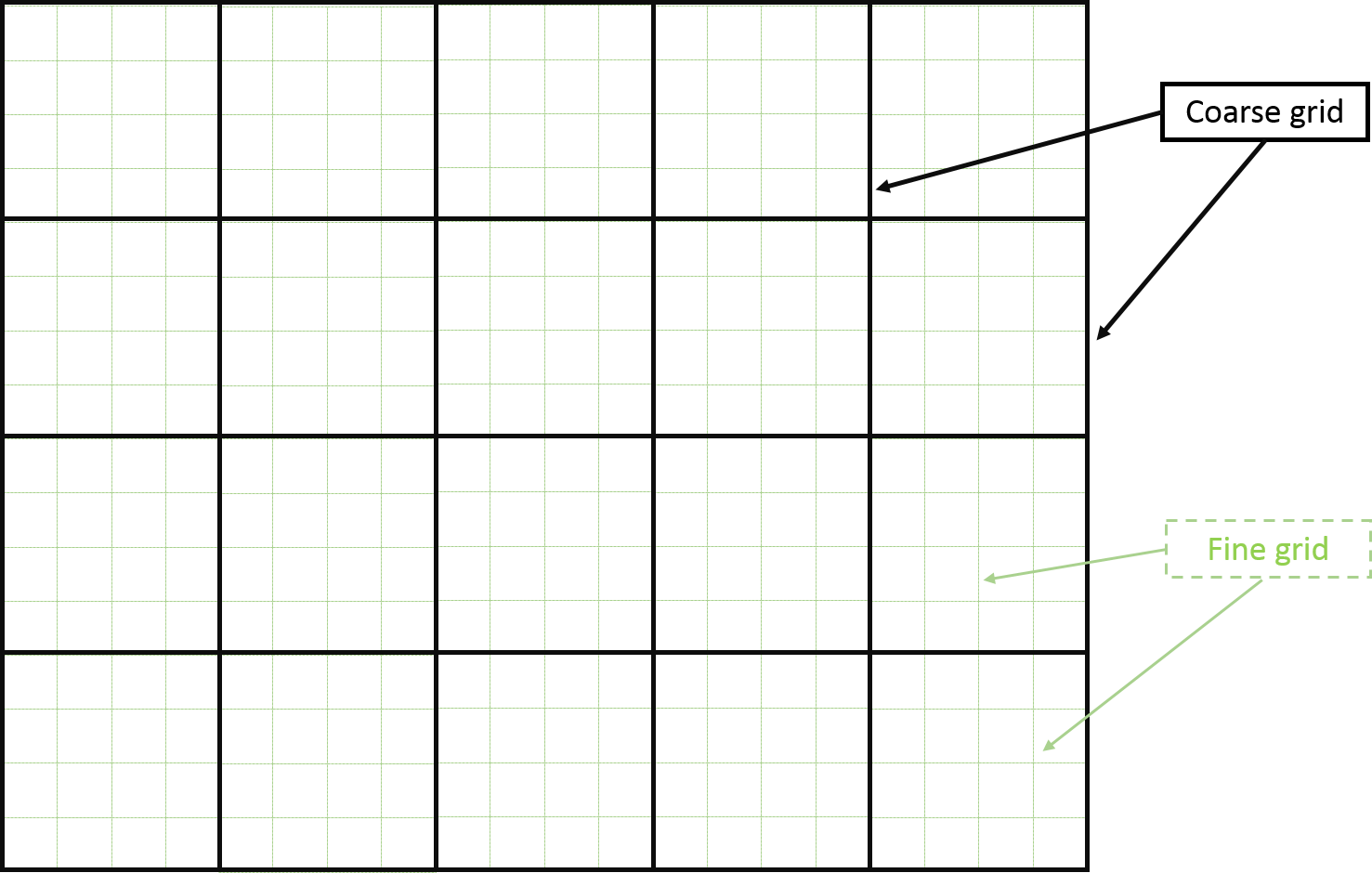}
			%\caption{fig1}
		\end{minipage}%
	}%
	\subfigure[coarse block and neighborhood.]{
		\begin{minipage}[t]{0.6\linewidth}
			\centering
			\includegraphics[width=3.8in]{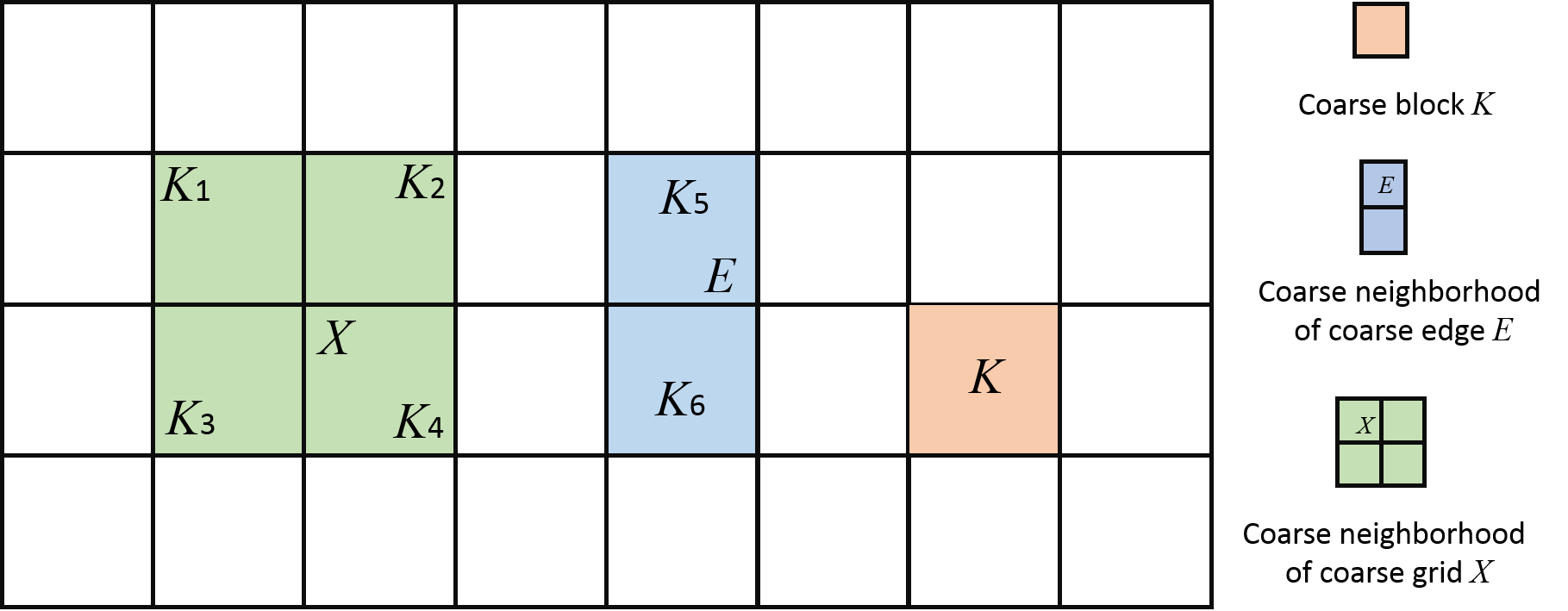}
			%\caption{fig2}
		\end{minipage}%
	}%
	\centering
	\caption{Illustration of mesh and neighborhood}
	\label{fig:neighborhood}
\end{figure}

In this section, we will derive the fine scale and mixed GMsFEM variational formulations for Problem \eqref{Pro:formula2}. 
Before introducing our method, we define the mesh and partitions needed in this paper. Let $\mathcal{T}^H$ be a standard conforming partition of the computational domain $\Omega$ into finite elements, where $H>0$ is the mesh size. We refer to this partition as the coarse-grid and assume that each coarse element is partitioned into a connected union of fine grid blocks. The fine grid will be denoted by $\mathcal{T}^h$, and is by definition a refinement of the coarse grid $\mathcal{T}^H$. We emphasize that we will use $K\in \mathcal{T}^H$ to denote a coarse element throughout the paper. Let 
$\mathcal{X}^H:=\{X_j\}_{j=1}^{N_c}$ be the set of nodes in the coarse grid $\mathcal{T}^H$, where $N_c$ is the number of the coarse nodes. Moreover, $\mathcal{X}^H_0$ is defined as a subset of $\mathcal{X}^H$ consisting of all interior coarse grid nodes.
We define the neighborhood $w_{X}$ of a coarse node $X\in \mathcal{X}_H$ by 
$$w_{X}: = \bigcup_j \{K_j\in \mathcal{T}^H|X\in \overline{K_j}\}.$$
Note that $w_{X}$ is the union of our all coarse elements $K_j\in \mathcal{T}^H$ sharing the coarse vertex $X$. Let $\mathcal{E}^H: = \{E_i\}_{i=1}^{N_e}$ be the set of all edges of the coarse mesh $\mathcal{T}^H.$ Furthermore, $\mathcal{E}^H_0$ is the subset of $\mathcal{E}^H$ containing all interior coarse edges. We define the coarse neighborhood $w_{E}$ of a coarse edge $E\in \mathcal{E}^H$ as the union of all coarse grid blocks having the edge $E$, namely,
$$w_{E} : = \bigcup_l \{K_l\in \mathcal{T}^H| E\in \partial K_l\}.$$ 
See Figure \ref{fig:neighborhood} for an illustration of neighborhoods of coarse edge and coarse grid. %We emphasize that we use $w_{X_i}$ or $w_{E_j}$ to denote a coarse neighborhood, and we use $K$ to denote a coarse element throughout the paper. We may use $w_i := w_{X_i}$ or $w_j := x_{E_j}$ for simplicity if no misunderstanding will be caused. 
For the time discretization, let $\Big\{T_j\Big\}_{j=0}^{j=J_t}$ 
$$0 = T_0<T_1<T_2<\cdots<T_{J_t} = T,$$
be a partition of $(0,T).$ In the following presentation, unknown with superscript $n$ equals its value at time $T_n$. For example, $p^{n} = p(\cdot, T_n). $

To introduce the variational formulation of Problem \eqref{Pro:formula2}, we define spaces $V, Z, Q, V^0$ and $Z^0$ as follows:
\begin{align*}
&	V = \bigg\{v\in \Big(H^1(\Omega)\Big)^d\bigg\},\quad Z=\bigg\{z\in H(div,\Omega)\bigg\}, \quad Q=\bigg\{q\in L^2(\Omega)\bigg\},\\
& V^0 = V\cap\Big\{v\in V|v = 0 \text{ on } (0,T]\times\partial\Omega\Big\}	,\quad Z^0 =Z\cap\Big\{z\in Z|z\cdot \vec{n} = 0 \text{ on }(0,T]\times\Gamma_2\Big\}.
\end{align*}
We first multiply \eqref{Pro:formula2a}, \eqref{Pro:formula2b} and \eqref{Pro:formula2c} with functions from $V^0$, $Z^0$ and $Q$, respectively. Next, applying Green's formula and making use of the boundary conditions on each portion, we get the variational formulation for Problem \eqref{Pro:formula2}: find $(u,g,p) \in (V,Z,Q)$ satisfying
\begin{subequations}
	\begin{align}
	&\int_{\Omega}-\nabla\cdot \sigma(u)v+\int_{\Omega}\alpha \nabla p v = 0,\forall v\in V^0,\label{Pro:formula5a}\\
	&\int_{\Omega}\kappa^{-1}\nu gz+\int_{\Omega}\nabla p z = 0,\forall z\in Z^0,\label{Pro:formula5b}\\
	&\int_{\Omega}\alpha \frac{\partial \nabla\cdot u}{\partial t}q+\int_{\Omega}\frac{1}{M}\frac{\partial p}{\partial t}q+\int_{\Omega}\nabla\cdot gq  = \int_{\Omega}fq, \forall q\in Q.\label{Pro:formula5c}
	\end{align}
	\label{Pro:formula5}
\end{subequations}

Let $V_h$ be the standard $\mathcal{Q}_1$ element for the approximation of fine-scale displacement $u$ on the fine grid $\mathcal{T}^h$, $Z_h$ be the standard lowest-order Raviart-Thomas space (RT$_0$) for fine-scale velocity $g$ approximation and $Q_h$ be the piecewise constant element for fine-scale pressure $p$ approximation for variational formulation \eqref{Pro:formula5}. Note that $V_h, Z_h$ and $Q_h$ are the fine-scale spaces and the corresponding solution set $(u_h,g_h,p_h)$ are used as our reference solutions in numerical experiments. $V_h^0$ and $Z_h^0$ can be defined similarly as $V^0$ and $Z^0$. Following same techniques as variational formulation of \eqref{Pro:formula5}, we have the fine-scale variational formulation: find $(u_h,g_h,p_h) \in (V_h,Z_h,Q_h)$ satisfying
\begin{subequations}
	\begin{align}
	&\int_{\Omega}-\nabla\cdot \sigma(u_h)v+\int_{\Omega}\alpha \nabla p_h v = 0,\forall v\in V_h^0,\label{Pro:formula3a}\\
	&\int_{\Omega}\kappa^{-1}\nu g_hz+\int_{\Omega}\nabla p_h z = 0,\forall z\in Z_h^0,\label{Pro:formula3b}\\
	&\int_{\Omega}\alpha \frac{\partial \nabla\cdot u_h}{\partial t}q+\int_{\Omega}\frac{1}{M}\frac{\partial p_h}{\partial t}q+\int_{\Omega}\nabla\cdot g_hq  = \int_{\Omega}fq, \forall q\in Q_h.\label{Pro:formula3c}
	\end{align}
	\label{Pro:formula3}
\end{subequations}

\begin{rem}
	Suppose $V_{\text{ms}}$, $Z_{\text{ms}}$, $Q_{\text{ms}}$ $V_{\text{ms}}^0$ and $Z_{\text{ms}}^0$ are some multiscale spaces for displacement, velocity and pressure which we will discuss later in Section \ref{Sec:multiscalebasis}. The variational formulation for our multiscale method is similar to fine scale formulation in \eqref{Pro:formula3c}. Therefore, the variational formulation for our multiscale method is: find $(u_{\text{ms}},g_{\text{ms}},p_{\text{ms}}) \in (V_{\text{ms}},Z_{\text{ms}},Q_{\text{ms}})$ satisfying	
	\begin{subequations}
	\begin{align}
	&\int_{\Omega}-\nabla\cdot \sigma(u_{\text{ms}})v+\int_{\Omega}\alpha \nabla p_{\text{ms}} v = 0,\forall v\in V_{\text{ms}}^0,\label{Pro:formula4a}\\
	&\int_{\Omega}\kappa^{-1}\nu g_{\text{ms}}z+\int_{\Omega}\nabla p_{\text{ms}} z = 0,\forall z\in Z_{\text{ms}}^0,\label{Pro:formula4b}\\
	&\int_{\Omega}\alpha \frac{\partial \nabla\cdot u_{\text{ms}}}{\partial t}q+\int_{\Omega}\frac{1}{M}\frac{\partial p_{\text{ms}}}{\partial t}q+\int_{\Omega}\nabla\cdot g_{\text{ms}}q  = \int_{\Omega}fq, \forall q\in Q_{\text{ms}}.	\label{Pro:formula4c}	
	\end{align}
	\label{Pro:formula4}
	\end{subequations}
\end{rem}
Finally, we notice that two terms $\alpha \frac{\partial \nabla\cdot u_h}{\partial t}$ and $\frac{1}{M}\frac{\partial p_h}{\partial t}$ in \eqref{Pro:formula3c} involve time derivative, which requires further discretization techniques. To facilitate our discussion, we define the following bilinear and linear operators:
\begin{align*}
&a(u,v) = \int_{\Omega} \sigma(u):\epsilon(v) dx,&& b(v,p) = \int_{\Omega}\alpha \nabla\cdot v p dx,\\
&c(q,v) = \int_{\Omega} \alpha q \nabla\cdot vdx,&&
d(q,p) = \int_{\Omega}\frac{1}{M}q pdx,\\
& e(q,g) = \int_{\Omega} q\nabla\cdot g dx,&& f(q) = \int_{\Omega} fq dx,\\
&j(z,g)=\int_{\Omega}\kappa^{-1}\nu zgdx,&& k(z,p)=\int_{\Omega}\nabla\cdot z p dx.
\end{align*}
One popular splitting method is fixed-stress splitting scheme.  The main idea is to combine \eqref{Pro:formula3b} and \eqref{Pro:formula3c} for the approximation of new step $p_{h}^{n+1}$ and $g_{h}^{n+1}$. Then pass the new $p_{h}^{n+1}$ to \eqref{Pro:formula3a} and calculate the new $u_{h}^{n+1}.$ In this way, \eqref{Pro:formula3a}-\eqref{Pro:formula3c} is divided as a sequence of variational formulations: find $(u_{h}^{n+1}, g_{h}^{n+1}, p_{h}^{n+1})\in V_h\times Z_h \times Q_h, n =0,1,2,\cdots, J_t-1,$ such that
\begin{subequations}
\begin{align}
a\Big(u_{h}^{n+1},v\Big)&=b(v,p_{h}^{n+1}),\forall v\in V_h^0,\label{splitting1:1}\\
j(z,g_{h}^{n+1})-k(z,p_{h}^{n+1})&=0,\forall z\in Z_h^0,\label{splitting1:2}\\
e(q,g_{h}^{n+1})+d(q,\frac{p_{h}^{n+1}}{\tau})&=f(q)-c(q,\frac{u_{h}^{n}-u_{h}^{n-1}}{\tau})+d(q,\frac{p_{h}^n}{\tau}),\forall q\in Q_h.\label{splitting1:3}
\end{align}
\label{splitting1}
\end{subequations}
In formulation \eqref{splitting1}, right hand sides of \eqref{splitting1:2}-\eqref{splitting1:3} only involve terms that can be computed at time step $n$ or before.  Only equations \eqref{splitting1:2}-\eqref{splitting1:3} are coupled in this scheme. We will use this splitting scheme in our numerical experiments in Section \ref{Sec:numerical results}.

Another feasible discretization method is the fully coupled method. All unknowns $u_{h}^{n+1}, g_{h}^{n+1}, h_{h}^{n+1}$ will be solved at a time in this method. Consequently, a much larger matrix will be created and it is time-consuming. The corresponding variational method is as follows: find $(u_{h}^{n+1}, g_{h}^{n+1}, p_{h}^{n+1})\in V_h\times Z_h \times Q_h, n =0,1,2,\cdots, J_t-1,$ such that
\begin{subequations}
\begin{align}
a(u_{h}^{n+1},v)-b(v,p_{h}^{n+1})&=0,\forall v\in V_h^0,\\
j(z,g_{h}^{n+1})-k(z,p_{h}^{n+1})&=0,\forall z\in Z_h^0,\\
c(q,\frac{u_{h}^{n+1}}{\tau})+e(q,g_{h}^{n+1})+d(q,\frac{p_{h}^{n+1}}{\tau})&=f(q)+c(q,\frac{u_{h}^{n}}{\tau})+d(q,\frac{p_{h}^n}{\tau}),\forall q\in Q_h.
\label{splitting2}
\end{align}
\end{subequations}
The variational formulation for multiscale space approximation of fixed-stress splitting and fully coupled splitting can be similarly derived.

\section{The Construction of Multiscale Basis Functions} \label{Sec:multiscalebasis}
As we have formed a sequence of variational formulation in Section \ref{Sec:variational and splitting}, we are left with the construction of multiscale spaces $V_{\text{ms}}$, $Z_{\text{ms}}$, $Q_{\text{ms}}$ $V_{\text{ms}}^0$ and $Z_{\text{ms}}^0$. The multiscale space $Q_{\text{ms}}$ is trivial in our method. It is piecewise constant with respect to the coarse partition $\mathcal{T}^H$. For the construction of multiscale space of velocity and displacement, the main idea is designing spectral problems to extract the dominant modes and thus get a reduced space. 
\begin{comment}
Before proceeding to the construction of $V_{\text{ms}}$ and $Z_{\text{ms}}$, we note that spaces $V_{\text{ms}}^0$ and $Z_{\text{ms}}^0$
are defined as follows:
$$	V_{\text{ms}}^0 = V_{\text{ms}}\cap\Big\{v\in V_{\text{ms}}|v = 0 \text{ on }\partial\Omega\Big\},\quad Z_{\text{ms}}^0 =Z_{\text{ms}}\cap\Big\{z\in Z_{\text{ms}}|z\cdot \vec{n} = 0\text{ on }\Gamma_2\Big\}.$$
%The critical point of multiscale method is reducing the number of basis by spectral problems. 
\end{comment}

\subsection{Multiscale space for velocity}
The multiscale space
$Z_{\text{ms}}$ is formed by solving a spectral problem on a snapshot space. The snapshot space of velocity are spanned by solutions of local problem with unit flux on part of local boundary.
 For an arbitrary coarse edge $E_i\in\mathcal{E}^H,$ suppose $E_i$ is the union of $l_i$ fine edges in $\mathcal{T}^h$, i.e., 
$E_i = \bigcup\limits_{j=1}^{l_i}e_j$, where $l_i$ is the total number of find-grid edges on $E_i$ and $e_j$ denotes a fine-grid edge in coarse edge $E_i$. For every fine edge $e_j\in E_i,$ we may define $l_i$ distinct fine edge delta functions on $E_i$ as follows:
\begin{equation*}
\delta^{j}_i =
\left\{
\begin{array}{lr}
1 \text{ on } e_j,&  \\
0 \text{ on } e_k, k\neq j.&
\end{array}
\right.
\end{equation*}
As indicated by the definition, $\delta^{j}_i$ is a piecewise constant function defined on $E_i$, and it has value $1$ on $e_j$ and $0$ on other fine-grid edges of $E_i$. Given these notations, we can define the following problem on the neighborhood $w_{E_i}$ of $E_i$: find $(g_i^j,p_i^j)\in (Z_h,Q_h)$ such that
\begin{equation}
\left\{
\begin{aligned}
\nabla p^j_i+\kappa^{-1}\nu g_i^j &=& 0&\, \text{ in } w_{E_i},\\
\nabla\cdot g_i^j&=& \alpha^j_i&\,\text{ in }w_{E_i},\\
g_i^j\cdot \vec{n}_i&=&0&\,\text{ on }\partial w_{E_i},\\
g_i^j\cdot \vec{m}_i&=& \delta^{j}_i&\,\text{ on }  E_i.
\label{pro:snapshotv}
\end{aligned}
\right.
\end{equation}
Here $\vec{n}_i$ denotes the outward unit normal vector on $\partial w_{E_i}$ and  $\vec{m}_i$ a fixed unit normal vector with respect to edge $E_i$.  $\alpha^j_i$ is yet to be determined. Indeed, the above problem can be solved separately on each coarse block of $w_{E_i}$. In this case, we construct $m_i$ corresponding to edge $E_i$.
  $\alpha^j_i$ is determined uniquely by the compatibility condition $\int_{K_l} \alpha^{(j)} = \int_{E} \delta^j_i,\quad \forall K_l \subset w_E.$
  
The collection of the solutions of the above local problems generates the snapshot space. We let $\Psi_j^{i,\text{snap}}: = g_i^{j}$ be the snapshot fields and define the snapshot space $Z_{\text{snap}}$ by
$$Z_{\text{snap}} = \text{ span }\{\Psi_j^{i,\text{snap}} : 1\leqslant j\leqslant l_i, 1\leqslant i\leqslant N_e\}.$$
\begin{comment}
To simplify notation, we will use the following single-index notation
$$Z_{\text{snap}} = \text{ span }\{\Psi_i^{\text{snap}} : 1\leqslant i\leqslant M_{\text{snap}}\},$$
where $M_{\text{snap}} = \sum_{i = 1}^{N_e} l_i$ is the total number of snapshot fields.
\end{comment} 
Moreover, we define the local snapshot space by 
$$Z_{\text{snap}}^i = \text{ span }\{\Psi_j^{i,\text{snap}} : 1\leqslant j\leqslant l_i\}.$$
Note that each $\Psi_i^{\text{snap}}$ is represented on the fine gird by the basis functions in $Z_h$. Therefore, each $\Psi_i^{\text{snap}}$ can be represented by a vector $\psi_i^{\text{snap}}$ containing the coefficients in the expansions of $\Psi_i^{\text{snap}}$ in the fine-grid basis functions. Then, we define
$$R_{\text{snap}} = [\psi_1^{\text{snap}},\dots,\psi_{M_{\text{snap}}}^{\text{snap}}],$$
which maps from the coarse space to the fine space.

\begin{comment}
\begin{rem}
We remark that the above local problem can be solved numerically on the underlying fine grid of $w_i$ by the lowest-order Raviart-Thomas element, so that the resulting velocity $v_j^{(i)} \in V_h$ (for simplicity, we keep the same notation for the discrete solution $v_j^{(i)}$).	
\end{rem}
\end{comment}

Next, we will perform a space reduction on the snapshot space through the use of some local spectral problems. The purpose of this is to determine the important modes in the snapshot space and to obtain a smaller space for approximating the solution. In the general setting, we consider the spectral problem of finding a real number $\lambda$ and a vector field $g\in Z_{\text{snap}}$ such that 
\begin{equation}
a(g,z) = \lambda s(g,z), \quad \forall z\in Z_{\text{snap}}
\label{pro:spectralv}
\end{equation}
where $a(g,z)$ and $s(g,z)$ are symmetric positive definite bilinear forms defined on $Z_{\text{snap}}\times Z_{\text{snap}}$. We consider $s(g,z)$ as an inner product on $Z_{\text{snap}}$ and define a linear operator $\mathcal{A}$: $Z_{\text{snap}}\rightarrow Z_{\text{snap}}$ by
$$s(\mathcal{A}g,z)=a(g,z).$$
The operator $\mathcal{A}$ has rapidly decaying eigenvalues if $\kappa$ is highly heterogeneous. Note that one can take $\mathcal{A}$ to be a compact operator. In practice, solving the above global spectral problem is inefficient. Therefore, the dimension reduction and the construction of the offline space are performed locally. In particular, the above spectral problem is solved for each $w_{E_i}$. An appropriate choice of spectral problem is vital for the final convergence. Below we list two effective and efficient spectral problems. In our numerical experiments, we will consider Spectral Problem 1.

\begin{itemize}
	\item \textbf{Spectral Problem 1}. We take 
	$$a_i(g,z) = \int_{E_i} \kappa^{-1} (g\cdot m_i)(z\cdot m_i),\quad s_i(g,z) = \int_{w_i}\kappa^{-1}g\cdot z + \int_{w_i}(\nabla \cdot g)(\nabla\cdot z).$$
	\item \textbf{Spectral Problem 2}. We take 
	$$a_i(g,z)=\int_{w_i}
	\kappa^{-1}g\cdot z,\quad s_i(g,z) = \int_{E_i} [p_g][p_z],$$
	where $(g,p_g)$ and $(z,p_z)$ are solutions of the local problem \eqref{pro:snapshotv}, and $[p]$ denotes the jump of the function $p$.
\end{itemize}
Without loss of generality, we assume the eigenpairs of spectral problem \eqref{pro:spectralv} can be sorted as 
$$\Big(\lambda_1^{(i)},\Theta_1^{(i)}\Big),\Big(\lambda_2^{(i)},\Theta_2^{(i)}\Big),\cdots,\Big(\lambda_{l_i}^{(i)},\Theta_{l_i}^{(i)}\Big),$$
$\Theta$
with $\Big(\lambda_k^{(i)}\Big)_{k=1}^{l_i}$ in a non-decreasing order. We will use the first $J_v^i$ eigenfunctions to form the offline space. The number $J_v^i$ depends on problem and will be chosen in the numerical experiments. Note that Using these eigenfunctions, offline basis functions can be constructed as
$$\Psi_k^{i,\text{off}} = \sum_{j=1}^{j=l_i} \Theta_{kj}^{(i)}\psi_j^{i,\text{snap}}, k = 1,2,\ddots,J_v^i.$$
The global offline space is then defined as
$$Z_{\text{ms}} = \text{span}\{\Psi_k^{i,\text{off}}:1\leqslant k\leqslant l_i,1\leqslant i\leqslant N_e\}.$$
To simplify notation, we will use the following single-index notation
$$Z_{\text{ms}} = \text{span}\{\Psi_k^{\text{off}}:1\leqslant k\leqslant M_{\text{off}}\},$$
where $M_{\text{ms}} = \sum_{i=1}^{N_e}l_i$ is the total number of offline basis functions. This space will be used as the approximation space for velocity in the GMsFEM system. Furthermore, we can define $Z_{\text{ms}}^0$ as the subspace of $Z_{\text{ms}}$ formed by the linear span of all the basis functions $\Psi_k^{\text{off}}$ corresponding to the set of all interior coarse edges $\mathcal{E}_H^0$. \\

\subsection{Multiscale space for displacement}
Similar to the construction of $Z_\text{ms}$ in previous subsection, we also select the most dominant modes to form the approximation space $V_\text{ms}$ by performing local spectral problems. One major difference however was the underlying spaces for spectral problems, where instead of constructing a snapshot space like $Z_{\text{snap}}$, we employ the space generated by all fine grid basis functions in $V_h$. To be specific, we want to find $(u,\lambda)\in \mathbb{R}\times V_h$ such that 
\begin{equation}
\hat{a}(u,v) = \lambda \hat{s}(u,v),\quad \forall v\in V_h.
\label{Pro:spectral2}
\end{equation}
Nevertheless, doing spectral problem on $V_h$ is memory-intensive and time-demanding. Therefore, we would prefer a reduced space $V_h(w_{X_j})$ to substitute $V_h$ in the spectral problem \eqref{Pro:spectral2}, where $V_h(w_{X_j})$ is the subspace of $V_h$ with domain restriction on a coarse neighborhood $w_{X_j}$. Hence, the spectral problem for displacement is written as: for every coarse neighborhood $w_{X_j}$, find $(u,\lambda)\in \mathbb{R}\times V_h(w_{X_j})$ such that 
\begin{equation}
\hat{a}(u,v) = \lambda \hat{s}(u,v),\quad \forall v\in V_h(w_{X_j}).
\end{equation}
In addition, the bilinear operators are chosen as 
 $$\hat{a}(v_m,v_n) = \int_{\Omega}\Big(2\mu\epsilon(v_m):\epsilon(v_n)+\lambda \nabla\cdot v_m \nabla\cdot v_n\Big),\quad \hat{s}(v_m,v_n) =\int_{\Omega}(\lambda+2\mu)v_m\cdot v_n, $$
$v_m, v_n\in V_h(w_{X_j}).$ Suppose $\Big\{(\mu_k^{(j)},\phi_k^{(j)})\Big\}_{k=1}^{k=l_j}$ are the eigenpairs of the problem \eqref{Pro:spectral2}, without loss of generality, we may assume they are arranged in a non-decreasing order by $\mu_k^{(j)}$. Suppose we intend to employ $J_u^j$ basis functions on $w_{X_j}$, then we can construct the corresponding offline basis functions as 
 $$\Phi_k^{j} = \sum_{m=1}^{J_u^j}v_m\phi_{km}^j, 1\leqslant j\leqslant N_c, 1\leqslant k\leqslant J_u^j,$$ 
 where $\phi_{km}^j$ are the $m$-th coordinates of $\phi_k^{(j)}$. According, we define local offline space $V_{\text{off}}^{w_j}$ as the space spanned by all the offline basis functions in the neighborhood $w_j$:
 $$
 V^{w_j}_{\text{off}}=: \text{ span } \{\Phi_k^{j}:1\leqslant k\leqslant J_u^i\}$$

To ensure the continuity of offline space, we multiply it by a multiscale partition of unity functions which are constructed by a local problem. For every coarse neighborhood of $w_{X_j}$, find $\xi^j_1 = \Big(\xi^j_{11},\xi^j_{12}\Big),$ $\xi^j_2 = (\xi^j_{21},\xi^j_{22})$ satisfying

\begin{equation*}
\left\{
 \begin{aligned}
\hat{a}(\xi^j_1,v) &=& 0 &\, \text{ in } w_{X_j},\\
\xi^j_{11} &= & g_j&\, \text{ on } \partial K, K\in w_{X_j},\\
\xi^j_{12} &= & 0&\, \text{ on } \partial K, K\in w_{X_j},\\
\xi^j_1 &= & 0&\, \text{ on } \partial w_{X_j},
\end{aligned}
\right.\quad
\left\{
\begin{aligned}
\hat{a}(\xi^j_2,v) &=& 0& \, \text{ in } w_{X_j},\\
\xi^j_{21}& = & 0&\, \text{ on } \partial K, K\in w_{X_j},\\
\xi^j_{22}& = & g_j&\, \text{ on } \partial K, K\in w_{X_j},\\
\xi^j_1 &= & 0&\, \text{ on } \partial w_{X_j}.
\end{aligned}
\right.
\end{equation*}
Here $g_j$ is a linear and continuous function on $\partial K$. The choice of $g_j$ can be referenced in paper \cite{hou1997multiscale}. For the sake of simplicity, we utilize the hat function regarding the coarse grid edge. Then our partition of unity multipliers are set as $\text{POU}_j = (\xi^j_{11},\xi^j_{22})$. Finally, we multiply the partition of unity functions by the eigenfunctions in the offline space $V_{\text{off}}^{w_j}$ to construct the resulting basis functions
 $$\Upsilon_{j,k} = \text{POU}_j\Phi_k^{j},\quad \text{for }1\leqslant j\leqslant N_c, 1\leqslant k\leqslant J_u^j.$$
Next, we define the multiscale space of displacement $V_{\text{ms}}$ as
$$V_{\text{ms}} = \text{span}\{\Upsilon_{j,k}: 1\leqslant j\leqslant N_c, 1\leqslant k\leqslant J_u^j\}.$$ Note that $V_{\text{ms}}^0$ are the subspace of $V_{\text{ms}}$ which excluded those basis generated by the coarse neighborhood on the boundary.
Once we constructed all the necessary multiscale spaces, we can use the splitting method and variational formulations introduced in Section \ref{Sec:variational and splitting} to get the final simulation.

\section{Numerical Results}\label{Sec:numerical results}
\begin{figure}[htbp]
	\centering
	\subfigure[$\kappa_1$.]{
		\begin{minipage}[t]{0.3\linewidth}
			\centering
			\includegraphics[width=2in]{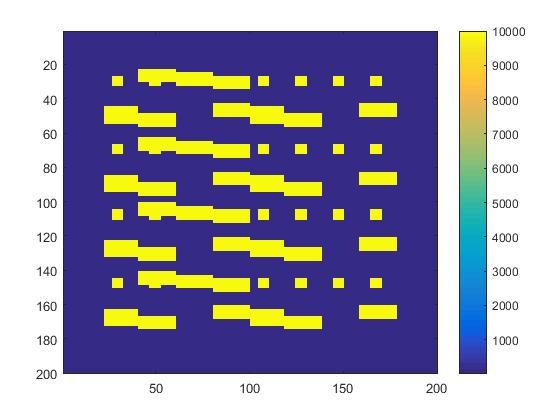}
			%\caption{fig1}
		\end{minipage}%
	}%
	\subfigure[$\kappa_2$.]{
		\begin{minipage}[t]{0.3\linewidth}
			\centering
			\includegraphics[width=2in]{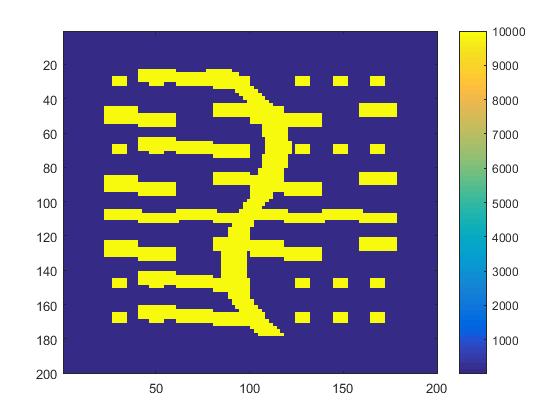}
			%\caption{fig2}
		\end{minipage}%
	}%
	\subfigure[$\kappa_3$.]{
	\begin{minipage}[t]{0.3\linewidth}
		\centering
		\includegraphics[width=2in]{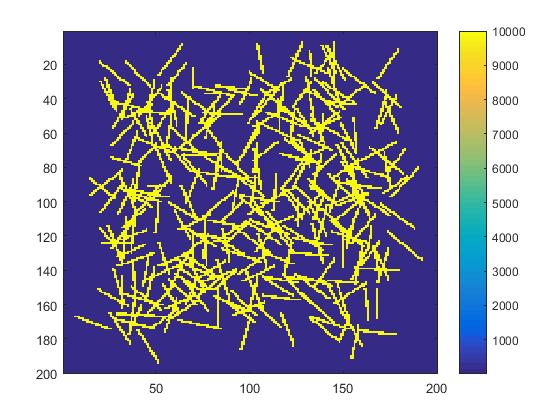}
		%\caption{fig2}
	\end{minipage}%
}%
	\centering
	\caption{Three high-contrast permeability fields used in the experiments}
	\label{fig:permeability}
\end{figure}

In this section, some numerical results are presented to illustrate the performance of our mixed GMsFEM for approximating problem \eqref{Pro:formula2}. In all simulations reported below, we employ the fixed-stress splitting scheme derived in Section \ref{Sec:variational and splitting}. The computational domain $D = (0,1)^2$. In our experiments, we will use three different permeability fields $\kappa_1$, $\kappa_2$ and $\kappa_3$.  Each permeability field can be divided into 2 subdomains based on heterogeneous coefficients. Figure \ref{fig:permeability} shows the subdomains distribution of $\kappa_1$, $\kappa_2$ and $\kappa_3$ used in our experiments. In Figure \ref{fig:permeability}, We choose $\kappa_i =1,i=1,2,3,$ in the blue region and $\kappa_i = 10^4,i=1,2,3,$ in the yellow region. Moreover, the coarse grid $\mathcal{T}^H$ and the fine grid $\mathcal{T}^h$ are $N\times N$ and $n\times n$ uniformly meshed, respectively. A fixed number of fine grid $n=200$ is employed, which means we may change the number of coarse grid but the number of total fine grid is set to be $200$. Other coefficients information are listed as follows.
\begin{enumerate}
	\item The Young's modulus $E$ is set to equal to the permeability field coefficient $\kappa$. 
	\item The Biot modulus $M$ equals $1$ in $\Omega_1$, and $10$ in $\Omega_2$.
	\item The Biot-Wills fluid-solid coupling coefficient $\alpha =0.9$. 
	\item The Poisson's ratio $\eta = 0.2$. 
	\item The Lam$\acute{e}$ coefficients $\lambda,\mu$ are determined by $\eta, E$ via relation \eqref{eqn:relationlambda}.
	\item The initial pressure $p_0(x,y) = xy(1-x)(1-y),\quad\forall (x,y)\in D.$
\end{enumerate}
Recall that we use a few multiscale basis functions on each coarse neighborhood $w_{X_i}$. These number of coarse basis determine the problem size (dimension of multiscale spaces, dim of $V_{\text{ms}},$). We assume that in each neighborhood, we select the same number of multiscale basis functions for velocity, i.e., $J_v^i = J_v$. Similarly, we choose equal number of basis functions for displacement, with $J_u^i = J_u$. Furthermore, 
we choose equal time step size $\tau$, i.e., $T_{i+1}-T_i=\tau, i=0,1\cdots J_t-1$. For simplicity of presentation, we introduce the following error quantities for displacement, velocity and pressure

\begin{align*}
&E_{L^2}^u = \frac{||u_{\text{ms}}(\cdot,T)-u_h(\cdot,T)||_{L^2\Omega}}{||u_h(\cdot,T)||_{L^2\Omega}},&&E_{a}^u = \frac{||u_{\text{ms}}(\cdot,T)-u_h(\cdot,T)||_{a}}{||u_h(\cdot,T)||_{a}}, \\
&E_{L^2}^g = \frac{||\frac{\kappa}{v}(v_{\text{ms}}(\cdot,T)-v_h(\cdot,T))||_{L^2\Omega}}{||\frac{\kappa}{\nu}v_h(\cdot,T)||_{L^2\Omega}}, &&E_{L^2}^p = \frac{||(p_{\text{ms}}(\cdot,T)-p_h(\cdot,T))||_{L^2\Omega}}{||p_h(\cdot,T)||_{L^2\Omega}},
\end{align*}
where $(u_{\text{ms}}(\cdot,T),v_{\text{ms}}(\cdot, T),p_{\text{ms}}(\cdot,T))$ is the multiscale solutions and $(u_{h}(\cdot,T),v_{h}(\cdot, T),p_{h}(\cdot,T))$ are the reference solution obtained by fine-scale solver. Note that $E_{L^2}^v$ is the weighted $L^2$ norm of velocity. 

\subsection{Model 1}
In the first simulation, we consider the proposed problem at $T=1$ with source term
$$ f(x,y,t)=\left\{
\begin{aligned}
2,\quad & x\in(0,\frac{1}{N}),y\in(0,\frac{1}{N}),t\in (0,1], \\
-2,\quad & x\in(\frac{N-1}{N},1),y\in(\frac{N-1}{N},1),t\in (0,1],
\end{aligned}
\right.
$$
and $\Gamma_1 = \emptyset$, i.e., 
$$g_0\cdot \vec{n} = 0\text{ on } \, (0,1]\times \partial D,\quad u_0 = 0 \text{ on } \, (0,1]\times \partial D.$$ In the following part, we call this explicit problem Model 1. We test Model 1 with permeability fields $\kappa=\kappa_1$ and $\kappa=\kappa_3$. Table \ref{tab:relationnumbasisu}-\ref{tab:relationnumoftimestep} exhibits the relationship between the approximation errors and number of displacement basis used in per coarse grid neighborhood $J_u$, number of time steps $J_t$ used and the number of velocity basis used in per coarse edge neighborhood $J_g$.
In Table \ref{tab:relationnumbasisu}, the errors $e_{L^2}^u$, $e_{H^1}^u$ of Model 1 at $T=1$ drop quickly if more displacement multiscale basis are employed. However, the convergence properties reach a plateau when enough displacement basis functions are used.
\begin{table}[htbp]
	\centering
	\caption{Convergence result for Model 1: Relationship between errors and $J_u$ with $\kappa_1$, $N=10,n=200, T=1, J_t =10, J_g=2$ }
	\label{tab:relationnumbasisu}
	\begin{tabular}{|r||r|r|r|r|}\hline
		$J_u$  & $e_{L^2}^u$ & $e_{H^1}^u$ & $e_{L^2}^p$ & $e_{L^2}^g$ \\  \hline
		4       & 0.3138 & 0.4862 & 0.0270 & 0.0801 \\
		8       & 0.0379 & 0.2534 & 0.0270 & 0.0801 \\
		12      & 0.0285 & 0.2365 & 0.0270 & 0.0801 \\
		16     & 0.0260 & 0.2303 & 0.0270 & 0.0801 \\
		20     & 0.0253 & 0.2267 & 0.0270 & 0.0801 \\
		24     & 0.0258 & 0.2240 & 0.0270 & 0.0801 \\\hline
	\end{tabular}
\end{table}
Meanwhile, we test the relationship of errors and number of velocity basis used per coarse neighborhood $J_g$. The result is shown in Table \ref{tab:relationnumbasisg}. We can see clearly that the error $e_{L^2}^g$ get smaller if we use more velocity basis, though at smaller scale. While other error estimators almost maintain the same level when $J_g$ changes. One possible reason is that the error $e_{L^2}^g$ is already small when $2$ multiscale basis of velocity is used in per coarse neighborhood.
% J_g
\begin{table}[htbp]
	\centering
		\caption{Convergence result for Model 1: Relationship between errors and $J_g$ with $\kappa_1$, $N=10,n=200, T=1, J_u =20, J_t=10$ }
	\label{tab:relationnumbasisg}
	\begin{tabular}{|r||r|r|r|r|}\hline
	$J_g$ & $e_{L^2}^u$ & $e_{H^1}^u$ & $e_{L^2}^p$ & $e_{L^2}^g$ \\  \hline
	2     & 0.0253 & 0.2267 & 0.0270 & 0.0801 \\
	3     & 0.0246 & 0.2261 & 0.0269 & 0.0573 \\
	4     & 0.0253 & 0.2258 & 0.0269 & 0.0377 \\
	5     & 0.0254 & 0.2257 & 0.0269 & 0.0304 \\
	6     & 0.0257 & 0.2257 & 0.0269 & 0.0245 \\
	\hline
\end{tabular}  	
\end{table}
%J_t
\begin{table}[htbp]
	\centering
	\caption{Convergence result for Model 1: Relationship between errors and $J_t$ with  $\kappa_1$, $N=10,n=200, T=1, J_u =20, J_g=2$}
	\label{tab:relationnumoftimestep}
	\begin{tabular}{|r||r|r|r|r|}\hline
		$J_t$  & $e_{L^2}^u$ & $e_{H^1}^u$ & $e_{L^2}^p$ & $e_{L^2}^g$ \\  \hline
		5    &0.0254	&0.2272&	0.0269&	0.0800\\
		10  &0.0253&	0.2267&	0.0270	&0.0801\\		
		20  &0.0253&	0.2266&	0.0270 &	0.0801\\
		40   &0.0253&	0.2266&	0.0270 &	0.0801\\\hline
	\end{tabular}
\end{table}
Simultaneously, we test several different time step sizes. The result is shown in Table \ref{tab:relationnumoftimestep}. The error quantities almost have no difference when we enlarge the number of time steps $J_t$. Similar results can be seen in Model 2 and thus we may fix $J_u, J_g, J_t$ as follows:$$J_u=20, J_g=2, J_t=10.$$
Table \ref{tab:model11} presents the results of $\kappa_1$ and $\kappa_3$. In both cases, the error  $e_{L^2}^u$, $e_{H^1}^u$,$e_{L^2}^p$ and  $e_{L^2}^v$ decrease rapidly when we enlarge the number of coarse grid. We see greater errors in displacement, velocity and pressure of $\kappa_3$ when compared with $\kappa_1$. The biggest possibility is that the heterogeneity properties in $\kappa_3$ is more complex.
\begin{table}
	\centering
	\caption{Convergence result of Model 1: Relationship between errors and $N$ with $n=200, T=1, J_u=20,J_t=10,J_v=2$}
	\begin{tabular}{|r||r|r|r|r||r|r|r|r|}
		\hline
	&	\multicolumn{4}{|c|}{$\kappa_1$} & \multicolumn{4}{|c|}{$\kappa_3$}\\
		\hline
				$N$  & $e_{L^2}^u$ & $e_{H^1}^u$ & $e_{L^2}^p$ & $e_{L^2}^v$& $e_{L^2}^u$ & $e_{H^1}^u$ & $e_{L^2}^p$ & $e_{L^2}^v$  \\\hline
		8      & 0.0303 & 0.2432 & 0.0478 & 0.1071 & 0.3732 & 0.5974 & 0.0383 & 0.3930\\
		10     & 0.0253 & 0.2267 & 0.0270 & 0.0801 & 0.1388 & 0.4539 & 0.0190 & 0.1215 \\
		20     & 0.0092 & 0.1456 & 0.0045 & 0.0496 & 0.0488 & 0.2844 & 0.0036 & 0.0655 \\
		25     & 0.0053 & 0.1222 & 0.0024 & 0.0356 & 0.0318 & 0.2398 & 0.0020 & 0.0542 \\\hline
		\hline
	\end{tabular}
\label{tab:model11}
\end{table}

\begin{figure}[H]
	\centering
	\subfigure[First component of $u$.]{
		\begin{minipage}[t]{0.3\linewidth}
			\centering
			\includegraphics[width=2in]{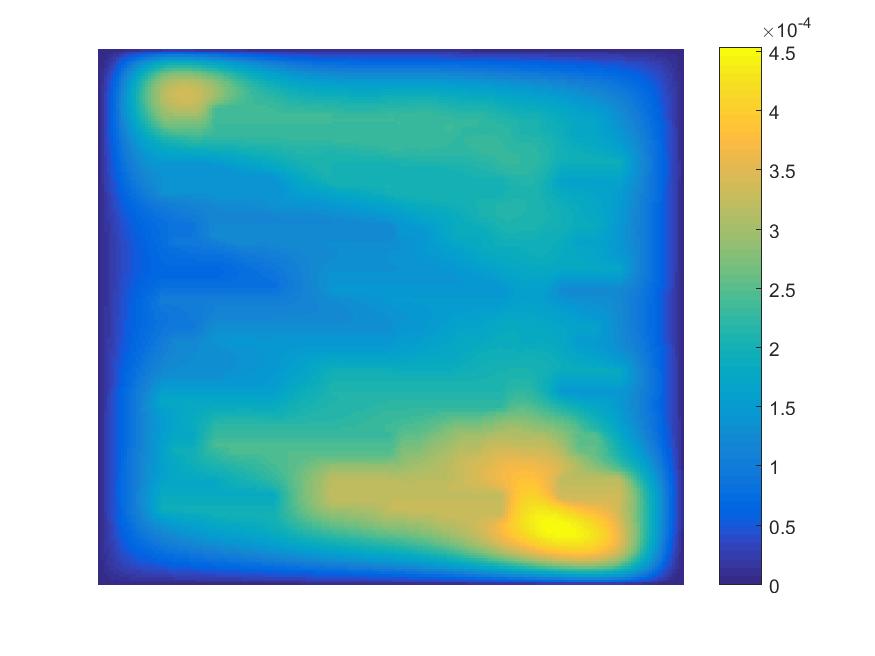}
			%\caption{fig1}
		\end{minipage}%
	}%
	\subfigure[Second component of $u$.]{
		\begin{minipage}[t]{0.3\linewidth}
			\centering
			\includegraphics[width=2in]{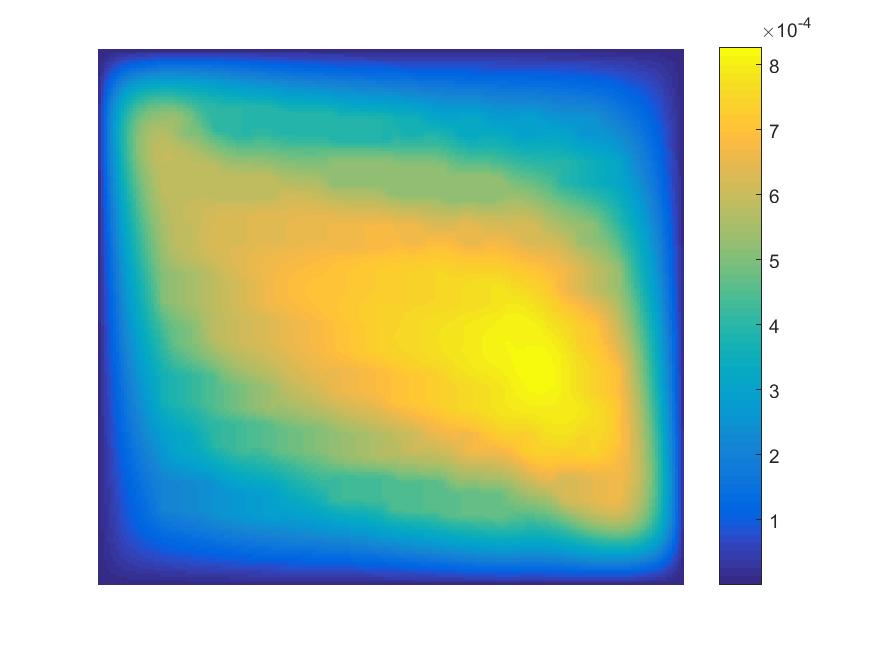}
			%\caption{fig2}
		\end{minipage}%
	}%
	\subfigure[Pressure $p$.]{
	\begin{minipage}[t]{0.3\linewidth}
		\centering
		\includegraphics[width=2in]{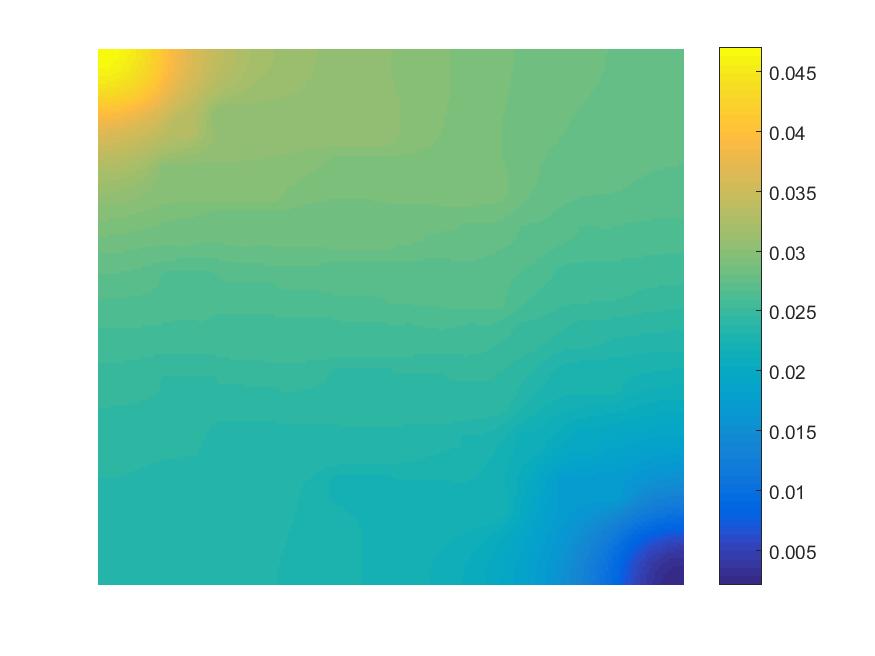}
		%\caption{fig2}
	\end{minipage}%
}%
	\centering
	\caption{Reference solution of Model 1 at $T=1$ with $\kappa = \kappa_1$}
	\label{image:reference1}
\end{figure}

\begin{figure}[H]
	\centering
	\subfigure[First component of $u$.]{
		\begin{minipage}[t]{0.3\linewidth}
			\centering
			\includegraphics[width=2in]{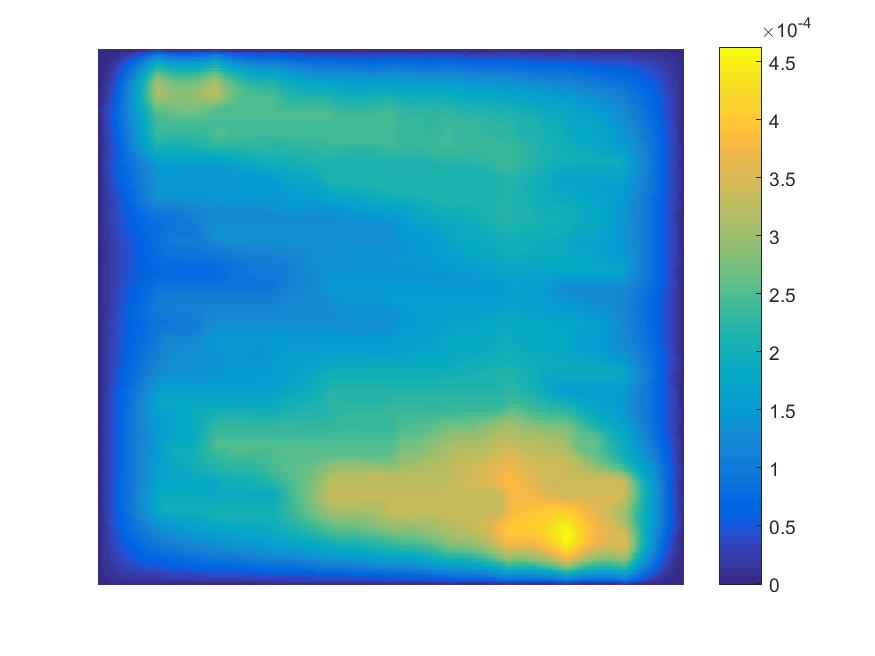}
			%\caption{fig1}
		\end{minipage}%
	}%
	\subfigure[Second component of $u$.]{
		\begin{minipage}[t]{0.3\linewidth}
			\centering
			\includegraphics[width=2in]{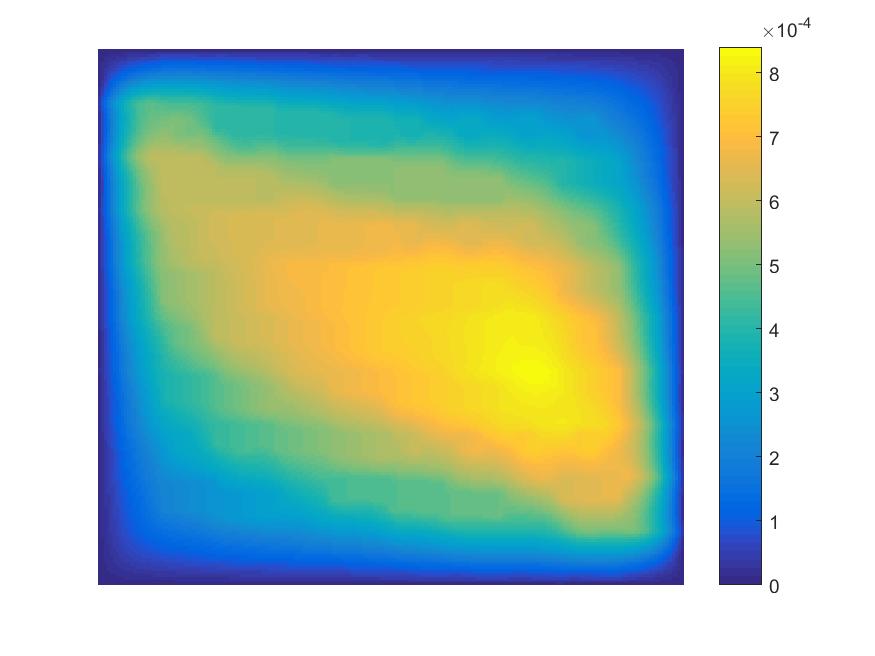}
			%\caption{fig2}
		\end{minipage}%
	}%
	\subfigure[Pressure $p$.]{
		\begin{minipage}[t]{0.3\linewidth}
			\centering
			\includegraphics[width=2in]{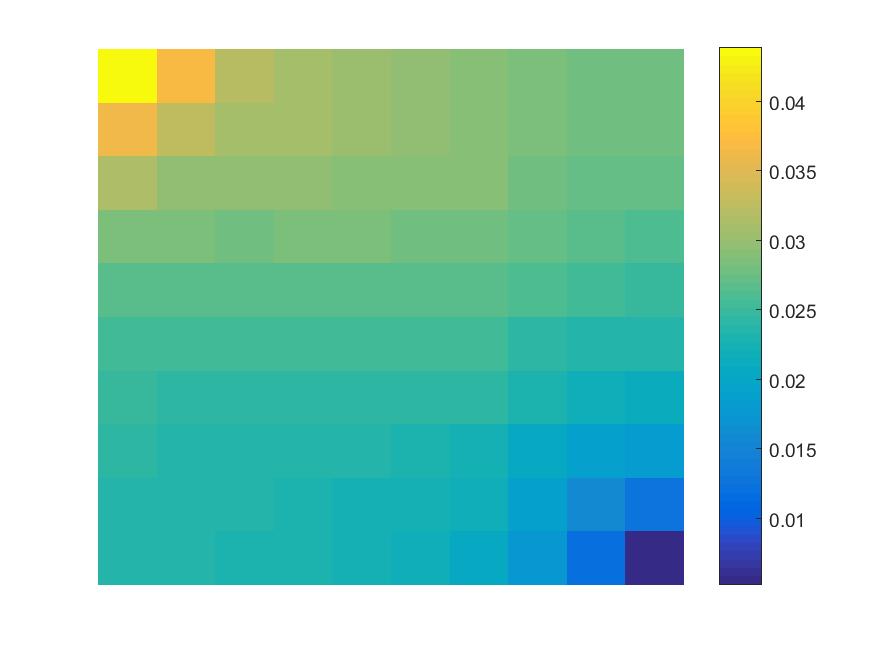}
			%\caption{fig2}
		\end{minipage}%
	}%
	\centering
	\caption{Mixed GMsFEM solution of Model 1 at $T=1$ with $\kappa = \kappa_1, N=10, n=200, J_u=20, J_g=2$}
	\label{image:mixed1}
\end{figure}

\begin{figure}[H]
	\centering
	\subfigure[Reference velocity $v$.]{
		\begin{minipage}[t]{0.3\linewidth}
			\centering
			\includegraphics[width=2in]{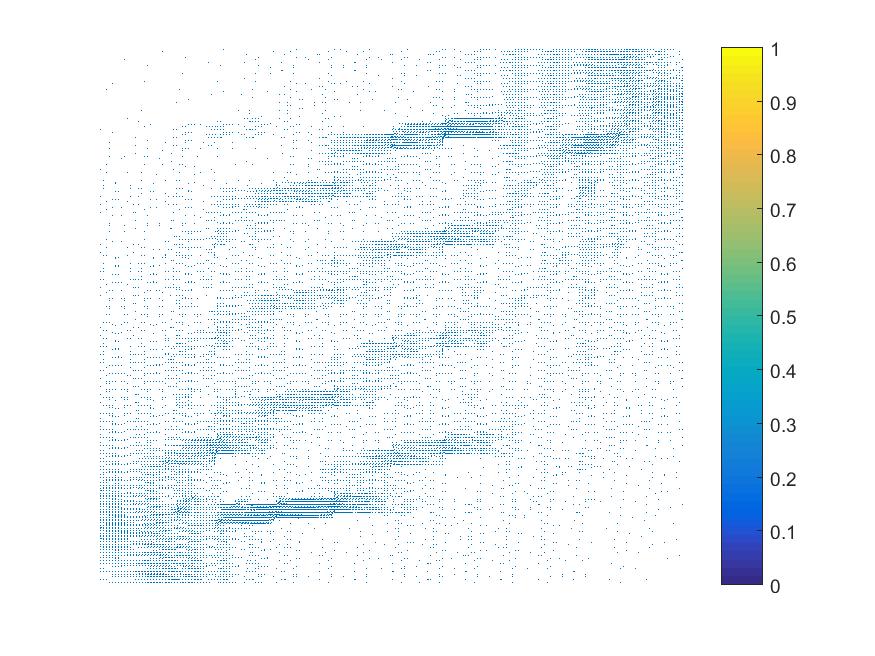}
			%\caption{fig1}
		\end{minipage}%
	}%
	\subfigure[Mixed velocity solution of $v$.]{
		\begin{minipage}[t]{0.3\linewidth}
			\centering
			\includegraphics[width=2in]{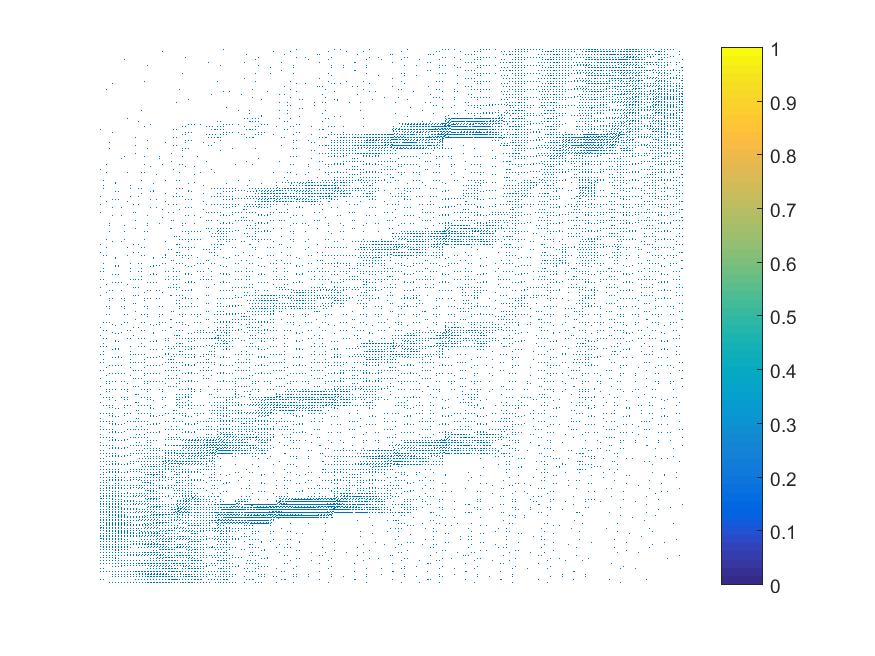}
			%\caption{fig2}
		\end{minipage}%
	}%
	\centering
	\caption{Comparison of reference solution and mixed GMsFEM solution for Model 1 at $T=1$ with $\kappa = \kappa_1$. Left: reference velocity solution. Right: mixed GMsFEM velocity solution with $N=10, n=200, J_u=20, J_g=2$}
	\label{image:reference11}
\end{figure}

Figure \ref{image:reference1}-\ref{image:reference11} are some images of our final result at $T=1$ of $\kappa_1$. Graphically, there is no observable difference between the reference solution and our mixed solution in this case.                         

\subsection{Model 2}
In the second model, we test the proposed method with $\Gamma_2 = \emptyset.$ We test with $T=1$. The boundary condition and source term are as follows:
$$f(x,y,t) = 1, \forall (x,y)\in D,t\in(0,1],\,u = 0,\,p = 0\text{ on } (0,1]\times \partial D.$$
 In the following part, we call this explicit problem Model 2. For Model 2, we will employ permeability fields $\kappa = \kappa_1$ and $\kappa= \kappa_2$. Relationships between the error quantities and $J_u$, $J_g$ and $J_t$ are similar to Model 1. Therefore, we choose the following numbers of basis:
$$J_u=20, J_g=2, J_t=10.$$
 Error results are shown in Table \ref{tab:model21}. For both $\kappa$, our scheme achieve good approximation. For $\kappa_1$, the $L^2$ error quantity for displacement dropped to $0.0164$ when there only 8 multiscale basis are chosen at each coarse neighborhood and the size of the coarse grid equals $\frac{1}{25}$. $\kappa_2$ are problem with more complex permeability media. Hence, the results are not as good as $\kappa_1$.
\begin{table}
	\centering
	\caption{Convergence result of Model 2: Relationship between errors and $N$ with $n=200,T=1, J_u=20,J_t=10,J_v=2$}
	\begin{tabular}{|r||r|r|r|r||r|r|r|r|}
		\hline
		&	\multicolumn{4}{|c|}{$\kappa_1$} & \multicolumn{4}{|c|}{$\kappa_2$}\\
		\hline
		$N$  & $e_{L^2}^u$ & $e_{H^1}^u$ & $e_{L^2}^p$ & $e_{L^2}^v$& $e_{L^2}^u$ & $e_{H^1}^u$ & $e_{L^2}^p$ & $e_{L^2}^v$  \\\hline
		8     & 0.3555 & 0.4869 & 0.2520 & 0.3882 & 0.4185 & 0.6161 & 0.2283 & 0.3045\\
		10    & 0.0985 & 0.3583 & 0.1438 & 0.0715 & 0.1856 & 0.4918 & 0.1564 & 0.0564\\
		20    & 0.0265 & 0.1872 & 0.0748 & 0.0463 & 0.0482 & 0.2522 & 0.0806 & 0.0347 \\
		25    & 0.0164 & 0.1503 & 0.0607 & 0.0332 & 0.0296 & 0.2011 & 0.0653 & 0.0311\\\hline
		\hline
	\end{tabular}
	\label{tab:model21}
\end{table}

\begin{figure}[H]
	\centering
	\subfigure[First component of $u$.]{
		\begin{minipage}[t]{0.3\linewidth}
			\centering
			\includegraphics[width=2in]{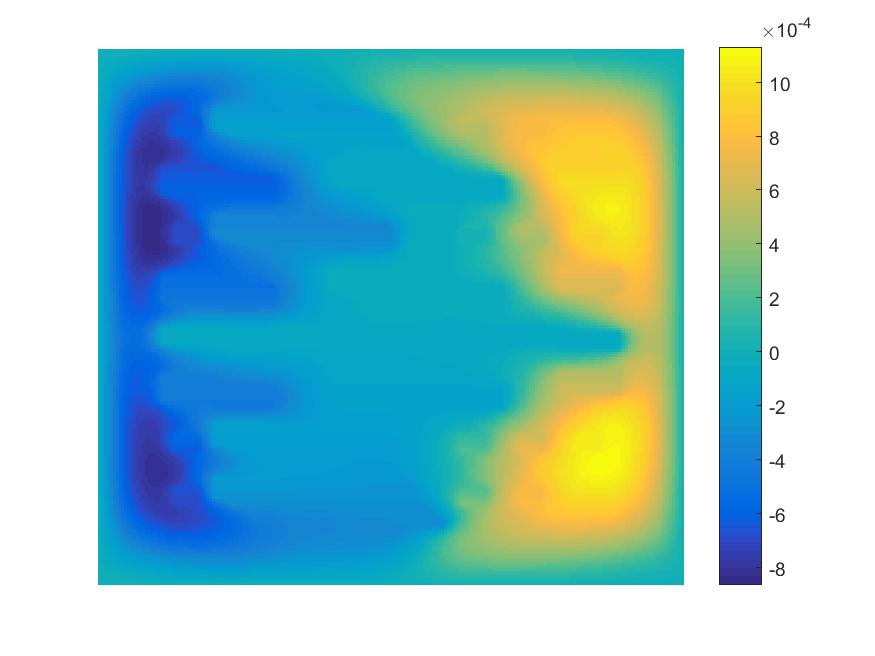}
			%\caption{fig1}
		\end{minipage}%
	}%
	\subfigure[Second component of $u$.]{
		\begin{minipage}[t]{0.3\linewidth}
			\centering
			\includegraphics[width=2in]{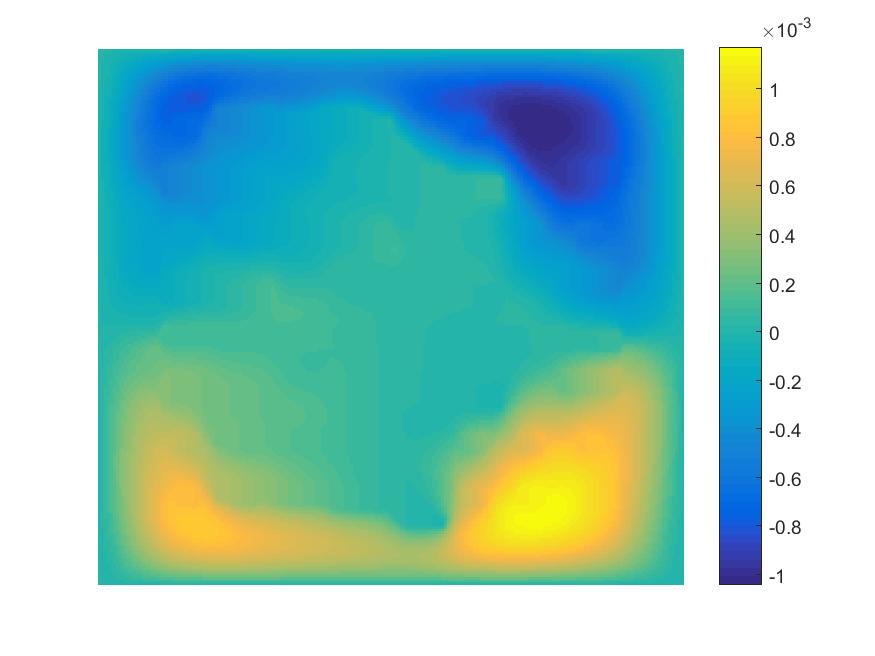}
			%\caption{fig2}
		\end{minipage}%
	}%
	\subfigure[Pressure $p$.]{
		\begin{minipage}[t]{0.3\linewidth}
			\centering
			\includegraphics[width=2in]{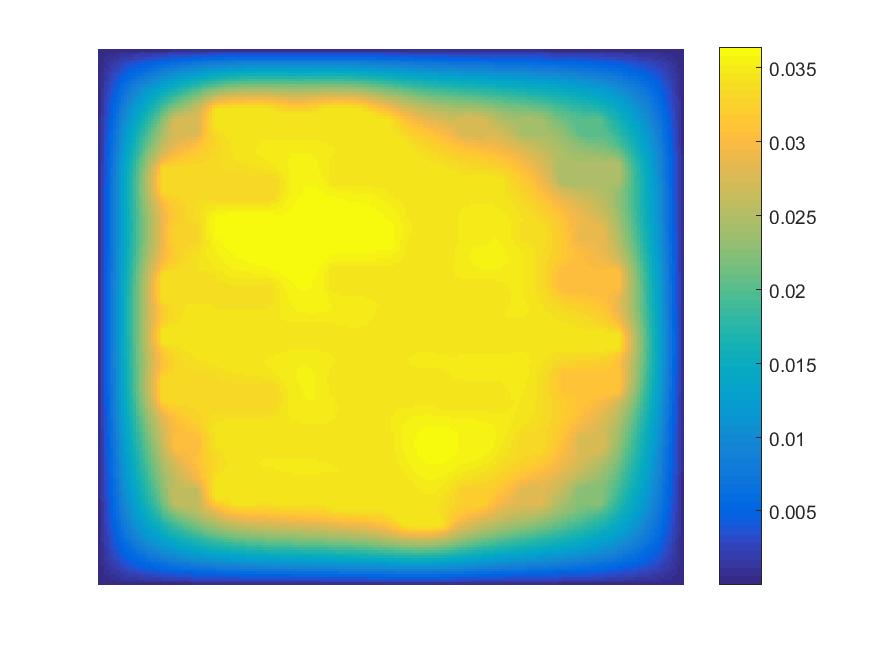}
			%\caption{fig2}
		\end{minipage}%
	}%
	\centering
	\caption{Reference solution of Model 2 at  $T=1$ with $\kappa = \kappa_2$}
	\label{image:reference2}
\end{figure}

\begin{figure}[H]
	\centering
	\subfigure[First component of $u$.]{
		\begin{minipage}[t]{0.3\linewidth}
			\centering
			\includegraphics[width=2in]{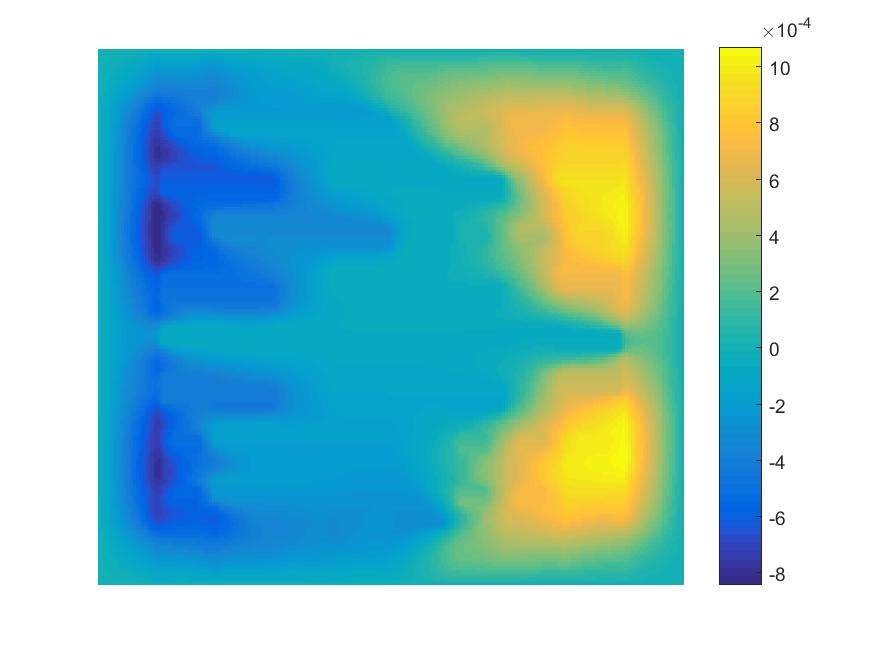}
			%\caption{fig1}
		\end{minipage}%
	}%
	\subfigure[Second component of $u$.]{
		\begin{minipage}[t]{0.3\linewidth}
			\centering
			\includegraphics[width=2in]{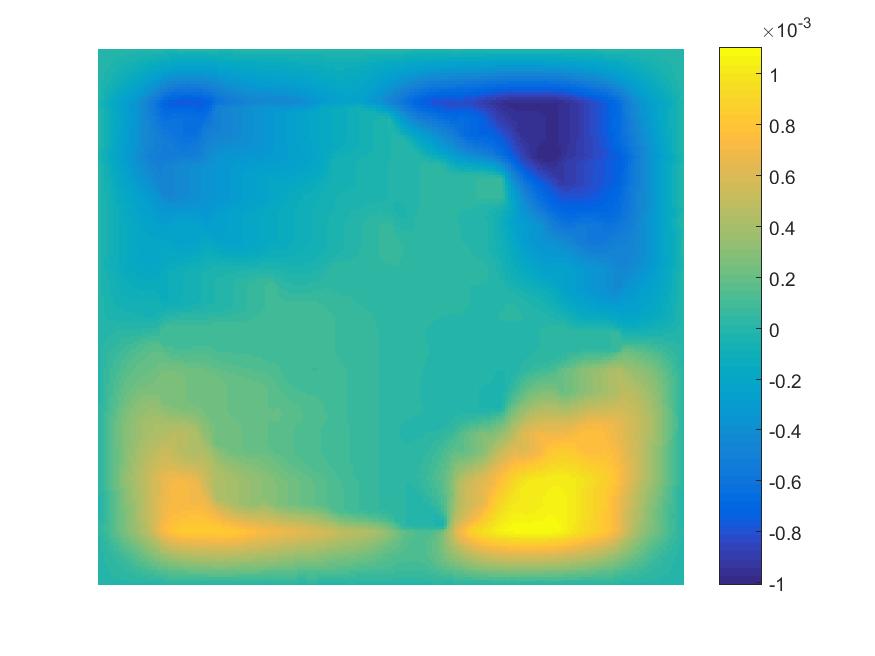}
			%\caption{fig2}
		\end{minipage}%
	}%
	\subfigure[Pressure $p$.]{
		\begin{minipage}[t]{0.3\linewidth}
			\centering
			\includegraphics[width=2in]{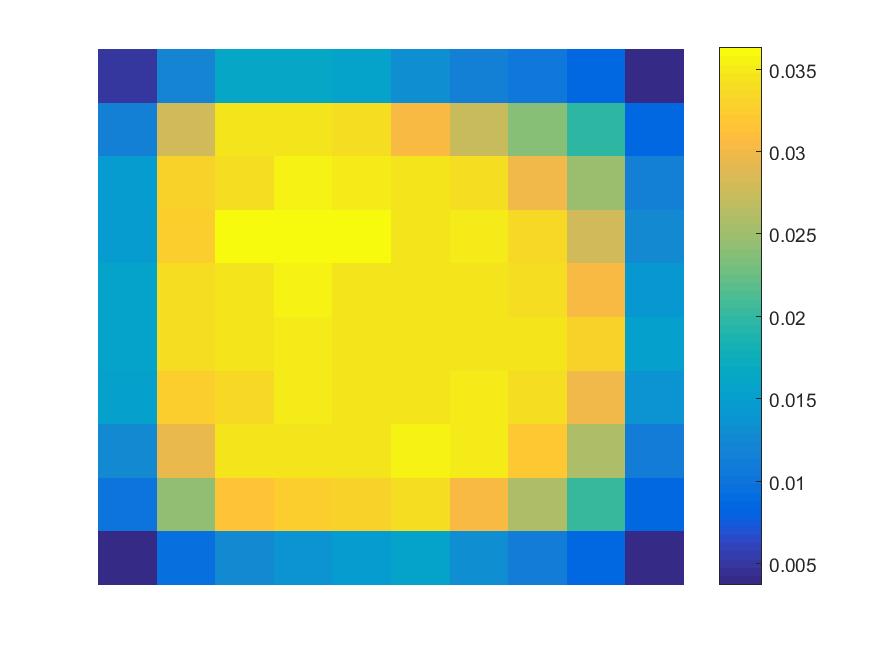}
			%\caption{fig2}
		\end{minipage}%
	}%
	\centering
	\caption{Mixed GMsFEM solution of Experiment 2 at $T=1$ with $\kappa = \kappa_2, N=10, n=200, J_u=20, J_g=2$}
	\label{image:mixed2}
\end{figure}

\begin{figure}[H]
	\centering
	\subfigure[Reference velocity $v$.]{
		\begin{minipage}[t]{0.3\linewidth}
			\centering
			\includegraphics[width=2in]{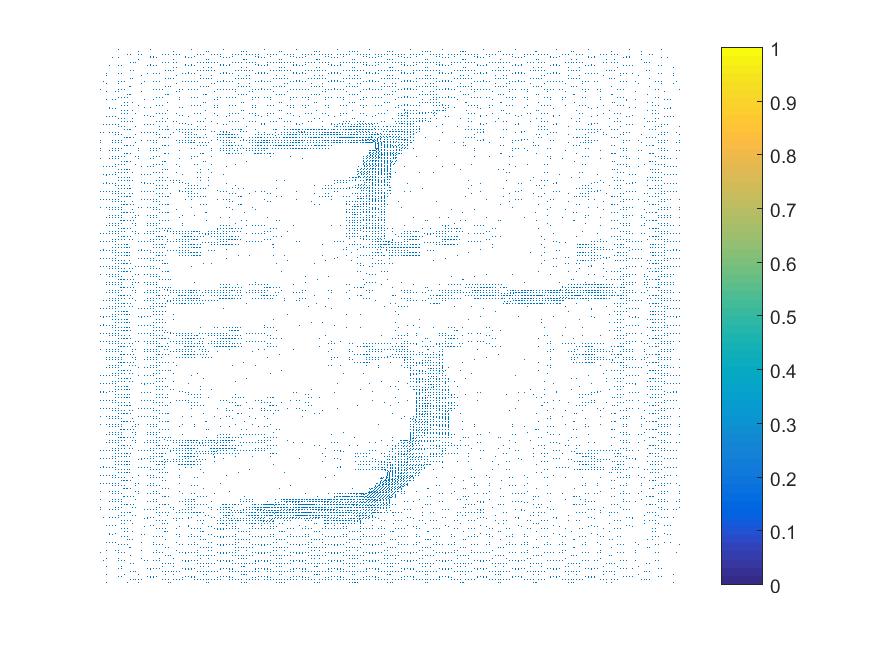}
			%\caption{fig1}
		\end{minipage}%
	}%
	\subfigure[Mixed velocity solution of $v$.]{
		\begin{minipage}[t]{0.3\linewidth}
			\centering
			\includegraphics[width=2in]{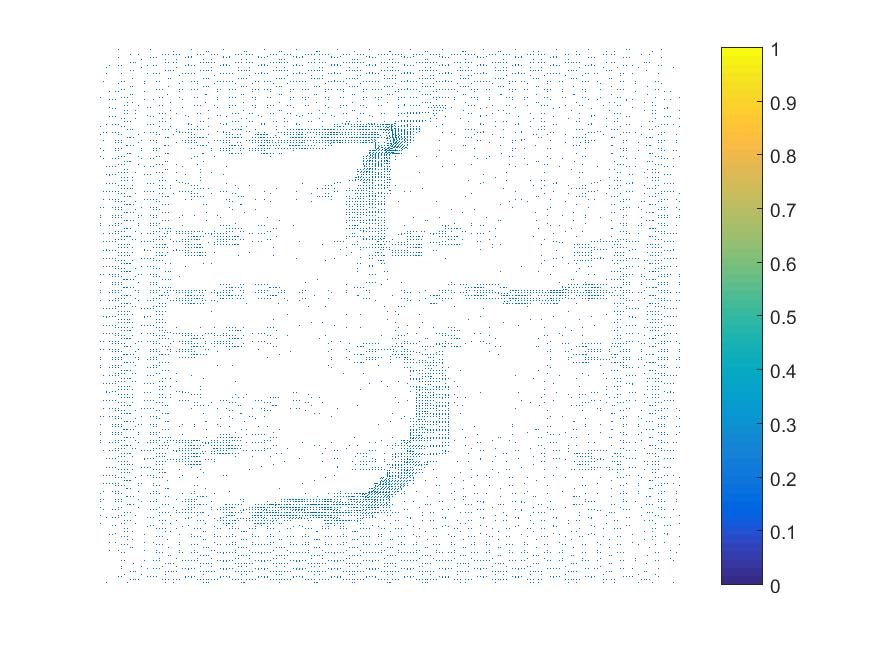}
			%\caption{fig2}
		\end{minipage}%
	}%
	\centering
	\caption{Comparison of reference solution and mixed GMsFEM solution for Model 2 at $T=1$ with $\kappa = \kappa_2$. Left: reference velocity solution. Right: mixed GMsFEM velocity solution with $ N=10, n=200, J_u=20, J_g=2$}
	\label{image:reference21}
\end{figure}

Figure \ref{image:reference2}-\ref{image:reference21} are some images of Model 2 with $\kappa = \kappa_2$. They demonstrate that our mixed GMsFEM works well on Model 2.

\section{Conclusion}\label{Sec:Conclusion}
In this paper, we have proposed a mass conservation method based on mixed finite element method and generalized multiscale finite element method (GMsFEM). We construct effective multiscale spaces by performing spectral problems for both velocity and displacement approximation. These multiscale basis functions are solutions of well designed local problems and can capture more heterogeneity properties of the medium. The numerical results show that our method works very well with only a few basis functions.
In the future, we will develop multiscale methods that are based on coupled basis functions for fluid velocity and elastic displacement.

\section*{Acknowledgement}

Eric Chung's work is partially supported by the Hong Kong RGC General Research Fund (Project numbers 14304719 and 14302018) and the CUHK Faculty of Science Direct Grant 2018-19.

\bibliographystyle{abbrv}
\bibliography{poroelasticity}

\begin{thebibliography}{10}

\bibitem{aarnes2004use}
J.~E. Aarnes.
\newblock On the use of a mixed multiscale finite element method for
  greaterflexibility and increased speed or improved accuracy in reservoir
  simulation.
\newblock {\em Multiscale Modeling \& Simulation}, 2(3):421--439, 2004.

\bibitem{abdulle2012heterogeneous}
A.~Abdulle, E.~Weinan, B.~Engquist, and E.~Vanden-Eijnden.
\newblock The heterogeneous multiscale method.
\newblock {\em Acta Numerica}, 21:1--87, 2012.

\bibitem{arbogast2007multiscale}
T.~Arbogast, G.~Pencheva, M.~F. Wheeler, and I.~Yotov.
\newblock A multiscale mortar mixed finite element method.
\newblock {\em Multiscale Modeling \& Simulation}, 6(1):319--346, 2007.

\bibitem{biot1941general}
M.~A. Biot.
\newblock General theory of three-dimensional consolidation.
\newblock {\em Journal of applied physics}, 12(2):155--164, 1941.

\bibitem{Biot1956}
M.~A. Biot.
\newblock Theory of propagation of elastic waves in a fluid-saturated porous
  solid. ii. higher frequency range.
\newblock {\em The Journal of the acoustical Society of america},
  28(2):179--191, 1956.

\bibitem{brown2016multiscale}
D.~L. Brown and D.~Peterseim.
\newblock A multiscale method for porous microstructures.
\newblock {\em Multiscale Modeling \& Simulation}, 14(3):1123--1152, 2016.

\bibitem{brown2016generalized}
D.~L. Brown and M.~Vasilyeva.
\newblock A generalized multiscale finite element method for poroelasticity
  problems i: linear problems.
\newblock {\em Journal of Computational and Applied Mathematics}, 294:372--388,
  2016.

\bibitem{chan2016adaptive}
H.~Y. Chan, E.~Chung, and Y.~Efendiev.
\newblock Adaptive mixed gmsfem for flows in heterogeneous media.
\newblock {\em Numerical Mathematics: Theory, Methods and Applications},
  9(4):497--527, 2016.

\bibitem{chen2003mixed}
Z.~Chen and T.~Hou.
\newblock A mixed multiscale finite element method for elliptic problems with
  oscillating coefficients.
\newblock {\em Mathematics of Computation}, 72(242):541--576, 2003.

\bibitem{chung2016adaptive}
E.~Chung, Y.~Efendiev, and T.~Y. Hou.
\newblock Adaptive multiscale model reduction with generalized multiscale
  finite element methods.
\newblock {\em Journal of Computational Physics}, 320:69--95, 2016.

\bibitem{chung2015mixed}
E.~T. Chung, Y.~Efendiev, and C.~S. Lee.
\newblock Mixed generalized multiscale finite element methods and applications.
\newblock {\em Multiscale Modeling \& Simulation}, 13(1):338--366, 2015.

\bibitem{chung2017onlineperforated}
E.~T. Chung, Y.~Efendiev, W.~T. Leung, M.~Vasilyeva, and Y.~Wang.
\newblock Online adaptive local multiscale model reduction for heterogeneous
  problems in perforated domains.
\newblock {\em Applicable Analysis}, 96(12):2002--2031, 2017.

\bibitem{chung2014adaptive}
E.~T. Chung, Y.~Efendiev, and G.~Li.
\newblock An adaptive gmsfem for high-contrast flow problems.
\newblock {\em Journal of Computational Physics}, 273:54--76, 2014.

\bibitem{chung2016mixedperforated}
E.~T. Chung, W.~T. Leung, and M.~Vasilyeva.
\newblock Mixed gmsfem for second order elliptic problem in perforated domains.
\newblock {\em Journal of Computational and Applied Mathematics}, 304:84--99,
  2016.

\bibitem{deng2017locally}
Q.~Deng, V.~Ginting, B.~McCaskill, and P.~Torsu.
\newblock A locally conservative stabilized continuous galerkin finite element
  method for two-phase flow in poroelastic subsurfaces.
\newblock {\em Journal of Computational Physics}, 347:78--98, 2017.

\bibitem{durlofsky1991numerical}
L.~J. Durlofsky.
\newblock Numerical calculation of equivalent grid block permeability tensors
  for heterogeneous porous media.
\newblock {\em Water resources research}, 27(5):699--708, 1991.

\bibitem{efendiev2013generalized}
Y.~Efendiev, J.~Galvis, and T.~Y. Hou.
\newblock Generalized multiscale finite element methods (gmsfem).
\newblock {\em Journal of Computational Physics}, 251:116--135, 2013.

\bibitem{efendiev2009multiscale}
Y.~Efendiev and T.~Y. Hou.
\newblock {\em Multiscale finite element methods: theory and applications},
  volume~4.
\newblock Springer Science \& Business Media, 2009.

\bibitem{ferronato2010fully}
M.~Ferronato, N.~Castelletto, and G.~Gambolati.
\newblock A fully coupled 3-d mixed finite element model of biot consolidation.
\newblock {\em Journal of Computational Physics}, 229(12):4813--4830, 2010.

\bibitem{fu2019computational}
S.~Fu, R.~Altmann, E.~T. Chung, R.~Maier, D.~Peterseim, and S.-M. Pun.
\newblock Computational multiscale methods for linear poroelasticity with high
  contrast.
\newblock {\em Journal of Computational Physics}, 395:286--297, 2019.

\bibitem{Gambolati2006}
G.~Gambolati, M.~Ferronato, and P.~Teatini.
\newblock Reservoir compaction and land subsidence.
\newblock {\em Revue européenne de génie civil}, 10:731--762, 09 2006.

\bibitem{gao2015numerical}
K.~Gao, E.~T. Chung, R.~L. Gibson~Jr, S.~Fu, and Y.~Efendiev.
\newblock A numerical homogenization method for heterogeneous, anisotropic
  elastic media based on multiscale theory.
\newblock {\em Geophysics}, 80(4):D385--D401, 2015.

\bibitem{hou1997multiscale}
T.~Y. Hou and X.-H. Wu.
\newblock A multiscale finite element method for elliptic problems in composite
  materials and porous media.
\newblock {\em Journal of computational physics}, 134(1):169--189, 1997.

\bibitem{kim2009stability}
J.~Kim, H.~A. Tchelepi, R.~Juanes, et~al.
\newblock Stability, accuracy and efficiency of sequential methods for coupled
  flow and geomechanics.
\newblock In {\em SPE reservoir simulation symposium}. Society of Petroleum
  Engineers, 2009.

\bibitem{kolesov2014splitting}
A.~E. Kolesov, P.~N. Vabishchevich, and M.~V. Vasilyeva.
\newblock Splitting schemes for poroelasticity and thermoelasticity problems.
\newblock {\em Computers \& Mathematics with Applications}, 67(12):2185--2198,
  2014.

\bibitem{Mura2016}
J.~Mura and A.~Caiazzo.
\newblock A two-scale homogenization approach for the estimation of porosity in
  elastic media.
\newblock In {\em Trends in Differential Equations and Applications}, pages
  89--105. Springer, 2016.

\bibitem{owhadi2007metric}
H.~Owhadi and L.~Zhang.
\newblock Metric-based upscaling.
\newblock {\em Communications on Pure and Applied Mathematics: A Journal Issued
  by the Courant Institute of Mathematical Sciences}, 60(5):675--723, 2007.

\bibitem{weinan2003heterognous}
E.~Weinan, B.~Engquist, et~al.
\newblock The heterognous multiscale methods.
\newblock {\em Communications in Mathematical Sciences}, 1(1):87--132, 2003.

\bibitem{wu2002analysis}
X.-H. Wu, Y.~Efendiev, and T.~Y. Hou.
\newblock Analysis of upscaling absolute permeability.
\newblock {\em Discrete and Continuous Dynamical Systems Series B},
  2(2):185--204, 2002.

\bibitem{yang2020online}
Y.~Yang, S.~Fu, and E.~T. Chung.
\newblock Online mixed multiscale finite element method with oversampling and
  its applications.
\newblock {\em Journal of Scientific Computing}, 82(2):31, 2020.

\bibitem{yang2019multiscale}
Y.~Yang, K.~Shi, and S.~Fu.
\newblock Multiscale hybridizable discontinuous galerkin method for flow
  simulations in highly heterogeneous media.
\newblock {\em Journal of Scientific Computing}, 81(3):1712--1731, 2019.

\bibitem{Zoback2010}
M.~D. Zoback.
\newblock {\em Reservoir geomechanics}.
\newblock Cambridge University Press, 2010.

\end{thebibliography}

\end{document}